\documentclass[a4paper,10pt,reqno]{article}
\usepackage{amsmath, amsfonts, amssymb, amsthm}

\input{epsf.sty}
\usepackage{mathrsfs}
\usepackage{amsmath}
\usepackage{amssymb}
\usepackage{graphicx}
\usepackage{subfigure}
\numberwithin{equation}{section}

\begin{document}

\title{On the Cauchy problem for a new integrable two-component system with peakon and weak kink solutions}

\author
{Kai $\mbox{Yan}^1$\footnote{E-mail: yankai419@163.com },\quad
Zhijun $\mbox{Qiao}^2$\footnote{E-mail: qiao@utpa.edu} \quad
and \quad Zhaoyang $\mbox{Yin}^1$\footnote{E-mail: mcsyzy@mail.sysu.edu.cn; zy409@cims.nyu.edu} \\
$^1\mbox{Department}$ of Mathematics,
Sun Yat-sen University,\\ Guangzhou, Guangdong 510275 , China \\
$^2\mbox{Department}$ of Mathematics, University of Texas-Pan American,\\ Edinburg, Texas 78541, USA}

\date{}
\maketitle

\begin{abstract}
This paper is contributed to study the Cauchy problem of a new integrable two-component
system with peaked soliton (peakon) and weak kink solutions. We first
establish the local well-posedness result for the Cauchy problem in Besov spaces, and then present a precise blow-up
scenario and a new blow-up result for strong solutions to the
system.\\


\end{abstract}

\section{Introduction}
In the past two decades, a large amount of literature was devoted to the celebrated Camassa-Holm
(CH) equation \cite{C-H}
\begin{eqnarray*}
m_t+u m_x+2u_x m+b u_x=0, \quad m=u-u_{xx},
\end{eqnarray*}
where $b$ is an arbitrary constant, which models the
unidirectional propagation of shallow water waves over a flat
bottom. Here $u(t, x)$ stands for the fluid velocity at time $t$
in the spatial $x$ direction \cite{C-H,D-G-H,J}. CH is also a
model for the propagation of axially symmetric waves in
hyperelastic rods \cite{Dai}. It has a bi-Hamiltonian structure
and is completely integrable \cite{C-H,F-F}. Also there is a
geometric interpretation of CH in terms of geodesic flow on the
diffeomorphism group of the circle \cite{C-K}. Its solitary waves
are peaked \cite{C-H-H,Qiao-CMP}. They are orbitally stable and
interact like solitons \cite{B-S-S,C-S1}. It is also worth
pointing out the peaked solitary waves replicate a feature that is
characteristic for the waves of great height-waves of largest
amplitude that are exact solutions of the governing equations for
water waves, cf. the discussion in the papers
\cite{Cinvent,C-Eann}.

The Cauchy problem and initial boundary value problem for CH have
been studied extensively \cite{C-Ep,C-Ec,D2,E-Y1}. It has been
shown that this equation is locally well-posed
\cite{C-Ep,C-Ec,D2} for initial data $u_0\in H^s(\mathbb{R})$,
$s>\frac{3}{2}$. Moreover, it has global strong solutions
\cite{Cf,C-Ep,C-Ec} and also finite time blow-up solutions
\cite{Cf,C-Ep,C-Ec,C-E}. On the other hand, it has global weak
solutions in $H^1(\mathbb{R})$ \cite{B-C2,C-M,X-Z}. The advantage
of CH in comparison with the KdV equation lies in the fact that CH
has peaked solitons and models wave breaking \cite{C-H-H,C-E} (by
wave breaking we understand that the wave remains bounded while
its slope becomes unbounded in finite time \cite{Wh}).

Another important integrable equation admitting peaked soliton (peakon) is the
famous Degasperis-Procesi (DP) equation \cite{D-P}
\begin{eqnarray*}
m_t+u m_x+3u_x m=0, \quad m=u-u_{xx},
\end{eqnarray*}
which is regarded as a model for nonlinear shallow water dynamics
\cite{CIL,C-L}. It was proved in \cite{D-H-H} that the DP equation
has a bi-Hamiltonian structure and an infinite number of
conservation laws, and admits exact peakon solutions which are
analogous to the CH peakons. The Cauchy problem and initial
boundary value problem for DP have been studied extensively in
\cite{Co-K,E-L-Y,E-L-Y2,E-Y1,L-Y,Y3,Y5}. Although DP is very
similar to CH in many aspects, especially in the structure of
equation, there are some essential differences between the two
equations. One of the famous features of DP equation is that it
has not only (periodic) peakon solutions \cite{D-H-H,Y5}, but also
(periodic) shock peakons \cite{E-L-Y2,Lu}. Besides, CH is a
re-expression of geodesic flow on the diffeomorphism group
\cite{C-K}, while the DP equation can be regarded as a non-metric
Euler equation \cite{E-K}.

Note that the nonlinear terms in the CH and DP equations are both quadratic. Indeed, there exist integrable
systems admitting peakons and nonlinear cubic terms, which are the modified CH equation
\begin{eqnarray}
m_t+((u^2-u^2_x)m)_x+b u_x=0, \quad m=u-u_{xx}
\end{eqnarray}
with a constant $b$, and the Novikov equation
\begin{eqnarray*}
m_t+u^2 m_x+3u u_x m=0, \quad m=u-u_{xx}.
\end{eqnarray*}
Eq.(1.1) was proposed independently in \cite{F-F,Fu,OR,Q1}. The
Lax pair and cusped soliton (cuspon) solutions for Eq.(1.1) have
been studied in \cite{Q1,Q2}. Recently, the formation of
singularities for Eq.(1.1), a wave-breaking mechanism and the
peakon stability of Eq.(1.1) in the case of $b=0$ were
investigated in \cite{Gui-CMP}. The Novikov equation was proposed
in \cite{NV1} and its Lax pair, bi-Hamiltonian structure, peakon
solutions, well-posedness, blow-up phenomena and global weak
solutions have been studied extensively in \cite{HH,HW,NV1,wu}.
Very recently, the following integrable equation with both
quadratic and cubic nonlinearity
\begin{eqnarray}
 m_t+\frac{1}{2}k_1((u^2-u^2_x)m)_x+\frac{1}{2}k_2( u m_x +2u_x m)+b u_x=0,
\end{eqnarray}
was investigated, where $m=u-u_{xx}$, $b$, $k_1$, and $k_2$ are
three arbitrary constants. Eq.(1.2) was first proposed by Fokas in
\cite{Fokas1}. Its Lax pair, bi-Hamiltonian structure, peakons,
weak kinks, kink-peakon interactional, and smooth soliton
solutions were studied recently in \cite{QXL}.

A natural idea is to extend such study to the multi-component generalized systems.
One of the most popular is the following integrable two-component Camassa-Holm shallow
water system
(2CH) \cite{C-L-Z,C-I}:
\begin{equation*}
\left\{\begin{array}{ll}
m_{t}+u m_x+2u_x m+\sigma\rho\rho_x=0,\\
\rho_{t}+(u\rho)_x=0,
\end{array}\right.
\end{equation*}
where $m=u-u_{xx}$ and $\sigma=\pm 1$, which becomes CH equation
when $\rho\equiv 0$. The Cauchy problems of (2CH) with $\sigma =
-1$ and with $\sigma = 1$ have been studied in \cite{E-Le-Y} and
\cite{C-I,G-Y,G-Y1,Gui-JFA}, respectively. Local well-posedness
for (2CH) with the initial data in Sobolev spaces and in Besov
spaces have been established in \cite{C-I,E-L-Y,Gui-JFA}. The
blow-up phenomena and global existence of strong solutions to
(2CH) in Sobolev spaces have been derived in
\cite{E-L-Y,G-Y,Gui-JFA}. The analyticity of solutions to the
Cauchy problem for (2CH) and  the initial boundary value problem
for (2CH) have been studied in \cite{Y-Y} and \cite{Y-Y3},
respectively. Recently, the existence of global weak solutions for
(2CH) with $\sigma = 1$ has been investigated in \cite{G-Y1}.
Another two famous systems are the modified two-component
Camassa-Holm system (M2CH) \cite{Holm}:
\begin{equation*}
\left\{\begin{array}{ll}
m_{t}+um_{x}+2u_x m +\sigma\rho\bar{\rho}_{x}=0,\\
\rho_{t}+(u\rho)_x=0,
\end{array}\right.
\end{equation*}
with $m=u-u_{xx}$, $\rho=(1-\partial^2_{x})(\bar{\rho}-\bar{\rho}_0)$, $\sigma =\pm 1$, and the
two-component Degasperis-Procesi system (2DP) \cite{Po2}:
\begin{equation*}
\left\{\begin{array}{ll}
m_t+u m_x+3u_x m+b\rho\rho_x=0,\\
\rho_{t}+u\rho_x+2u_x\rho=0,
\end{array}\right.
\end{equation*}
with $m=u-u_{xx}$ and a real constant $b$. When $\rho\equiv 0$,
(M2CH) and (2DP) are reduced to CH and DP equations, respectively.
The Cauchy problems and initial boundary value problems for (M2CH)
and (2DP) have been studied in many works, for example
\cite{G1,G2,T3,Y-Y2,Y-Y3} and \cite{Y-Y1,Y-Y3}, respectively.
It is noted to point out that the nonlinear terms in above three
two-component systems are all quadratic.

In this paper, we consider the following integrable two-component system with both quadratic and cubic nonlinearity
proposed in \cite{XQ}:
\begin{equation}
\left\{\begin{array}{ll}
m_t+\frac{1}{2}[(uv-u_x v_x)m]_x-\frac{1}{2}(uv_x-u_x v)m+bu_x=0,\\
n_t+\frac{1}{2}[(uv-u_x v_x)n]_x+\frac{1}{2}(uv_x-u_x v)n+bv_x=0,
\end{array}\right.
\end{equation}
where $m=u-u_{xx}$, $n=v-v_{xx}$, and $b$ takes an arbitrary value. System (1.3) is reduced to the CH
equation, the modified cubic CH equation Eq.(1.1), and the generalized CH equation Eq.(1.2) as $v=2$, $v=2u$,
and $v=k_1u+k_2$, respectively. Integrability of this system, its bi-Hamiltonian structure, and
infinitely many conservation laws were already presented in \cite{XQ}. In addition, this system is the first two-component system admitting weak kink solutions.
Let us set up the Cauchy problem for the above system as follows:
\begin{equation}
\left\{\begin{array}{ll}
m_t+\frac{1}{2}(uv-u_x v_x)m_x=-\frac{1}{2}(u_x n+v_x m)m+\frac{1}{2}(uv_x-u_x v)m-bu_x,\\
n_t+\frac{1}{2}(uv-u_x v_x)n_x=-\frac{1}{2}(u_x n+v_x m)n-\frac{1}{2}(uv_x-u_x v)n-bv_x,\\
m(0,x)=m_0(x), \\
n(0,x)=n_0(x).
\end{array}\right. \label{mn}
\end{equation}
By using the same approach in \cite{Y-Y}, the analytic solutions to System (\ref{mn}) can be readily proved in
both variables, globally in space and locally in time. However, the goal of this paper is to establish the local well-posedness regime for the Cauchy problem in Besov spaces, present the precise blow-up scenario and a new blow-up result for strong solutions to the system, and derive the peakon and weak kink solutions.

With regard to the locally well-posed problem, we first apply the classical Kato semigroup theory to obtain the local well-posedness for System (\ref{mn}) with initial data $(m_0, n_0)$ belonging in Sobolev space $H^s(\mathbb{R})\times H^s(\mathbb{R})$ as $s\geq 1$. Subsequently, we take advantage of the transport equation theory, Littlewood-Paley's decomposition and some fine estimates of Besov spaces to establish the local well-posedness for System (\ref{mn}) in Besov spaces (see Theorem 3.2 below), which particularly implies the system is locally well-posed with initial data $(m_0, n_0)\in H^s(\mathbb{R})\times H^s(\mathbb{R})$ for $\frac 1 2<s\neq \frac 3 2$. This almost improves the corresponding result by using Kato's semigroup approach.

In order to investigate the blow-up phenomena, we here make good use of the fine structure of System (\ref{mn}). It is not difficult to verify that System (\ref{mn}) possesses the following two conservation laws:
$$ H_1= \frac{1}{2} \int_{\mathbb{R}} (u v+u_x v_x) d x
=\frac{1}{2} \int_{\mathbb{R}} u n d x
=\frac{1}{2} \int_{\mathbb{R}} v m d x,$$
$$ H_2= \frac{1}{4} \int_{\mathbb{R}} \left((u^2 v_x+u^2_x v_x-2u u_x v)n+2 b(u v_x-u_x v)\right) d x.$$
Indeed, as mentioned before, this system has infinitely many conservation laws.
But, unfortunately, one cannot seem to seek and then utilize some of them to control the quantities
$||u(t,\cdot)||_{L^\infty}$ and $||v(t,\cdot)||_{L^\infty}$
directly. While it is very important to bound them in studying the blow-up phenomena of System (\ref{mn}).
This difficulty has been overcome in some
sense by exploiting the characteristic ODE related to System (\ref{mn})
to establish some invariant properties of the solutions and making good use of the structure of the system
itself (we here need to treat two cross terms $\frac{1}{2}(uv_x-u_x v)m$ and $-\frac{1}{2}(uv_x-u_x v)n$
in the system), see Lemma 4.2 below. It is noted to point out that, for the special case Eq.(1.1) as $b=0$,
with the conserved quantity $\int_{\mathbb{R}} u m d x=||u||^2_{H^1(\mathbb{R})}$ in hand,
one can directly control $||u(t,\cdot)||_{L^\infty}$ by using Sobolev's embedding theorem, that is,
$||u(t,\cdot)||_{L^\infty(\mathbb{R})}\leq C ||u(t,\cdot)||_{H^1(\mathbb{R})}
= C ||u_0(\cdot)||_{H^1(\mathbb{R})}$, for all $t\in[0,T)$.
And thus the blow-up phenomena of Eq.(1.1) with $b=0$ has been studied in \cite{Gui-CMP}.
On the other hand, in view of the uselessness of the
conservation laws of System (\ref{mn}) again (mainly because the regularity is not high
enough), we directly investigate the transport equation in terms
of the flow $\frac{1}{2}(u_x n+v_x m)$ which is the slope of
$\frac{1}{2}(u v-u_x v_x)$ (see Lemma 4.3 below) to derive a new
blow-up result with respect to initial data (see Theorem 4.3
below). Overall, we haven't used any conservation laws rather than the almost symmetrical structure of System (\ref{mn}) in the whole paper.

The rest of our paper is organized as follows. In Section 2, we will recall some facts on the Littlewood-Paley decomposition, the nonhomogeneous Besov spaces and their some useful properties, and the transport equation theory as well. In Section 3, we establish the local well-posedness of System (\ref{mn}). In Section 4, we derive the precise blow-up scenario and a new blow-up result of strong solutions to System (\ref{mn}). Section 5 is devoted to discussing peakon and weak kink solutions.

\section{Preliminaries}
\newtheorem {remark2}{Remark}[section]
\newtheorem{theorem2}{Theorem}[section]
\newtheorem{lemma2}{Lemma}[section]
\newtheorem{definition2}{Definition}[section]
\newtheorem{proposition2}{Proposition}[section]

In this section, we will recall some facts on the Littlewood-Paley decomposition, the nonhomogeneous Besov spaces and their some useful properties, and the transport equation theory, which will be used in the sequel.

\begin{proposition2}
\cite{BCD} (Littlewood-Paley decomposition) There exists a couple of smooth functions $(\chi,\varphi)$ valued in $[0,1]$, such that $\chi$ is supported in the ball $B\triangleq \{\xi\in\mathbb{R}^n:|\xi|\leq \frac 4 3\}$, and $\varphi$ is supported in the ring $C\triangleq \{\xi\in\mathbb{R}^n:\frac 3 4\leq|\xi|\leq \frac 8 3\}$. Moreover,
$$\forall\,\ \xi\in\mathbb{R}^n,\,\ \chi(\xi)+{\sum\limits_{q\in\mathbb{N}}\varphi(2^{-q}\xi)}=1,$$
and
$$\textrm{supp}\,\ \varphi(2^{-q}\cdot)\cap \textrm{supp}\,\ \varphi(2^{-q'}\cdot)=\emptyset,\,\ \text{if}\,\ |q-q'|\geq 2,$$
$$\textrm{supp}\,\ \chi(\cdot)\cap \textrm{supp}\,\ \varphi(2^{-q}\cdot)=\emptyset,\,\ \text{if}\,\ q\geq 1.$$
\end{proposition2}

Then for all $u \in \mathcal{S'}$, we can define the nonhomogeneous dyadic blocks as follows. Let
$$\Delta_q{u}\triangleq 0,\,\ if\,\ q\leq -2,$$
$$\Delta_{-1}{u}\triangleq \chi(D)u=\mathcal{F}^{-1}\chi \mathcal{F}u,$$
$$\Delta_q{u}\triangleq \varphi(2^{-q}D)u=\mathcal{F}^{-1}\varphi(2^{-q}\xi)\mathcal{F}u,\,\ if \,\ q\geq 0.$$
Hence,
$$ u={\sum\limits_{q\in\mathbb{Z}}}\Delta_q{u} \quad  in \,\  \mathcal{S'}(\mathbb{R}^n),$$
where the right-hand side is called the nonhomogeneous Littlewood-Paley decomposition of $u$.

\begin{remark2}
(1) The low frequency cut-off $S_q$ is defined by
$$S_q{u}\triangleq {\sum\limits_{p=-1}^{q-1}}\Delta_p{u}=\chi(2^{-q}D)u=\mathcal{F}^{-1}\chi(2^{-q}\xi)\mathcal{F}u,
\,\ \forall\,\ q\in\mathbb{N}.$$
(2) The Littlewood-Paley decomposition is quasi-orthogonal in $L^2$ in the following sense:
$$\Delta_p\Delta_q{u}\equiv 0,\quad \text{if} \,\ |p-q|\geq2,$$
$$\Delta_q(S_{p-1}u \Delta_p{v})\equiv 0,\quad \text{if} \,\ |p-q|\geq5,$$
for all $u,v \in \mathcal{S'}(\mathbb{R}^n)$.\\
(3) Thanks to Young's inequality, we get
$$||\Delta_q{u}||_{L^p},\,\ ||S_q{u}||_{L^p}\leq C||u||_{L^p},\quad \forall\,\ 1\leq p\leq\infty,$$
where $C$ is a positive constant independent of $q$.
\end{remark2}

\begin{definition2}
\cite{BCD} (Besov spaces) Let $s\in\mathbb{R}, 1\leq p,r\leq\infty$. The nonhomogeneous Besov space $B^s_{p,r}(\mathbb{R}^n)$ ($B^s_{p,r}$ for short) is defined by
$$B^s_{p,r}(\mathbb{R}^n)\triangleq \{f \in \mathcal{S'}(\mathbb{R}^n):||f||_{B^s_{p,r}}< \infty\},$$
where
$$||f||_{B^s_{p,r}}\triangleq ||2^{qs}\Delta_q{f}||_{l^r(L^p)}=||(2^{qs}||\Delta_q{f}||_{L^p})_{q\geq-1}||_{l^r}.$$
If $s=\infty$, $B^\infty_{p,r}\triangleq \bigcap\limits_{s\in\mathbb{R}}B^s_{p,r}$.
\end{definition2}

\begin{definition2}
Let $T>0$, $s\in\mathbb{R}$ and $1\leq p \leq\infty$. Set
$$E^s_{p,r}(T)\triangleq C([0,T]; B^s_{p,r})\cap C^1{([0,T]; B^{s-1}_{p,r})},\quad \text{if}\,\ r<\infty,$$
$$E^s_{p,\infty}(T)\triangleq L^{\infty}([0,T]; B^s_{p,\infty})\cap Lip\,([0,T]; B^{s-1}_{p,\infty})$$
and
$$E^s_{p,r}\triangleq \bigcap\limits_{T>0}E^s_{p,r}(T).$$
\end{definition2}

\begin{remark2}
By Definition 2.1 and Remark 2.1(3), we can deduce that
$$||\Delta_q{u}||_{B^s_{p,r}},\,\ ||S_q{u}||_{B^s_{p,r}}\leq C||u||_{B^s_{p,r}},$$
where $C$ is a positive constant independent of $q$.
\end{remark2}

In the following proposition, we list some important properties of Besov spaces.
\begin{proposition2}
\cite{BCD,Gui-JFA} Suppose that $s\in\mathbb{R}, 1\leq p,r,p_{i},r_{i}\leq\infty, i=1,2.$  We have\\
(1) Topological properties: $B^s_{p,r}$ is a Banach space which is continuously embedded in $\mathcal{S'}$.\\
(2) Density: $\mathcal{C}^{\infty}_c$ is dense in $B^s_{p,r}$ $\Longleftrightarrow 1\leq p,r < \infty.$\\
(3) Embedding: $B^s_{p_1,r_1}\hookrightarrow B^{s-n({1\over p_1}-{1\over p_2})}_{p_2,r_2}$, \,\  if\quad $p_1\leq p_2$ and  $r_1\leq r_2$,
$$B^{s_2}_{p,r_2}\hookrightarrow B^{s_1}_{p,r_1}\quad  \text{locally compact},\quad \text{if}\,\ \, s_1 < s_2.$$
(4) Algebraic properties: $\forall s>0$, $B^s_{p,r}\bigcap L^\infty$ is an algebra.
Moreover, $B^s_{p,r}$ is an algebra, provided that $s>{n\over p}$ or $s\geq{n\over p}\,\ and\,\ r=1$.\\
(5) 1-D Morse-type estimates:\\
(i) $\forall s_{1}\leq \frac{1}{p} < s_2$ ($s_{2}\geq {1\over p}$ if $r=1$) and $s_1+s_2>0$, we have
\begin{eqnarray}
||fg||_{B^{s_1}_{p,r}(\mathbb{R})}\leq C||f||_{B^{s_1}_{p,r}(\mathbb{R})}||g||_{B^{s_2}_{p,r}(\mathbb{R})}.
\end{eqnarray}
(ii)For $s>0$,
\begin{eqnarray}||fg||_{B^s_{p,r}(\mathbb{R})}\leq C(||f||_{B^s_{p,r}(\mathbb{R})}||g||_{L^\infty(\mathbb{R})}+||g||_{B^s_{p,r}(\mathbb{R})}||f||_{L^\infty(\mathbb{R})}).
\end{eqnarray}
(iii) In Sobolev spaces $H^s(\mathbb{R})= B^s_{2,2}(\mathbb{R})$, we have for $s>0$,
\begin{eqnarray}||f\partial_x g||_{H^s(\mathbb{R})}\leq
C(||f||_{H^{s+1}(\mathbb{R})}||g||_{L^\infty(\mathbb{R})}+||f||_{L^\infty(\mathbb{R})} ||\partial_x g||_{H^s(\mathbb{R})}),
\end{eqnarray}
where $C$ is a positive constant independent of $f$ and $g$.\\
(6) Complex interpolation:
\begin{eqnarray}
\quad \quad \quad ||f||_{B^{\theta {s_1}+(1-\theta){s_2}}_{p,r}}\leq ||f||^{\theta}_{B^{s_1}_{p,r}}||f||^{1-\theta}_{B^{s_2}_{p,r}}, \quad \forall u\in B^{s_1}_{p,r}\cap B^{s_1}_{p,r},\quad \forall \theta \in[0,1].
\end{eqnarray}
(7) Fatou lemma: if $(u_n)_{n\in \mathbb{N}}$ is bounded in $B^s_{p,r}$ and $u_n \to u $ in $\mathcal{S'}$, then $u\in B^s_{p,r}$ and
$$||u||_{B^s_{p,r}}\leq \liminf\limits_{n\to \infty} ||u_n||_{B^s_{p,r}}.$$
(8) Let $m\in\mathbb{R}$ and $f$ be an $S^m$-multiplier (i.e., $f: \mathbb{R}^n \to \mathbb{R}$ is smooth and satisfies that for all $\alpha\in \mathbb{N}^n$, there exists a constant $C_\alpha$ such that $|\partial^\alpha{f(\xi)}|\leq C_\alpha(1+|\xi|)^{m-|\alpha|}$ for all $\xi \in \mathbb{R}^n$). Then the operator $f(D)$ is continuous from $B^s_{p,r}$ to $B^{s-m}_{p,r}.$
\end{proposition2}

Now we state some useful results in the transport equation theory, which are crucial to the proofs of our main theorems later.
\begin{lemma2}
\cite{BCD} (A priori estimates in Besov spaces) Let $1\leq p,r\leq \infty$ and $s>-\min ({1\over p}, 1-{1\over p}).$ Assume that $f_0\in B^s_{p,r}$, $F\in L^1(0,T; B^s_{p,r})$, and $\partial_x v$ belongs to $L^1(0,T; B^{s-1}_{p,r})$
if $s> 1+{1\over p}$ or to $L^1(0,T; B^{1\over p}_{p,r}\cap L^\infty)$ otherwise.
If $f\in L^\infty(0,T; B^s_{p,r})\bigcap C([0,T]; \mathcal{S'})$ solves the following 1-D linear transport equation:
\[(T)\left\{\begin{array}{l}
\partial_t f+v\,\partial_x f=F,\\
f|_{t=0} =f_0,
\end{array}\right.\]
then there exists a constant $C$ depending only on $s,p$ and $r$, and such that the following statements hold:\\
(1) If $r=1$ or $s\neq 1+{1\over p}$,
\begin{equation*}
||f(t)||_{B^s_{p,r}}\leq ||f_0||_{B^s_{p,r}}\,+\, \int_0^t ||F(\tau)||_{B^s_{p,r}}d\tau\,+\, C\int_0^t V'(\tau)||f(\tau)||_{B^s_{p,r}} d\tau
\end{equation*}
or hence,
\begin{equation*}
||f(t)||_{B^s_{p,r}}\leq e^{CV(t)} (||f_0||_{B^s_{p,r}}\,+\, \int_0^t e^{-CV(\tau)} ||F(\tau)||_{B^s_{p,r}}d\tau)
\end{equation*}
with $V(t)=\int_0^t ||\partial_x v(\tau)||_{B^{1\over p}_{p,r}\cap L^\infty}d\tau$ if $s<1+{1\over p}$ and $V(t)=\int_0^t ||\partial_x v(\tau)||_{B^{s-1}_{p,r}}d\tau$ else.\\
(2) If $v=f$, then for all $s>0$, (1) holds true with $V(t)=\int_0^t ||\partial_x v(\tau)||_{L^\infty}d\tau $.\\
(3) If $r<\infty$, then $f\in C([0,T]; B^s_{p,r})$. If $r=\infty$, then $f\in C([0,T]; B^{s'}_{p,1})$ for all $s'<s$.
\end{lemma2}

\begin{lemma2}
\cite{Gui-JFA} (A priori estimate in Sobolev spaces) Let $0<\sigma<1$. Assume that
$f_0\in H^\sigma$, $F\in L^1(0,T; H^\sigma)$, and $v,\partial_x v\in L^1(0,T; L^\infty)$.
If $f\in L^\infty(0,T; H^\sigma)\bigcap C([0,T]; \mathcal{S'})$ solves $(T)$,
then $f\in C([0,T]; H^{\sigma})$, and there exists a constant $C$ depending only on $\sigma$ such that the following statement holds:
\begin{equation*}
||f(t)||_{H^{\sigma}}\leq ||f_0||_{H^{\sigma}}\,+\, \ C \int_0^t ||F(\tau)||_{H^{\sigma}}d\tau\,+\, C\int_0^t V'(\tau)||f(\tau)||_{H^{\sigma}} d\tau
\end{equation*}
or hence,
\begin{equation*}
||f(t)||_{H^{\sigma}}\leq e^{CV(t)} (||f_0||_{H^{\sigma}}\,+\, \int_0^t ||F(\tau)||_{H^{\sigma}}d\tau)
\end{equation*}
with $V(t)=\int_0^t (||v(\tau)||_{L^\infty}+||\partial_x v(\tau)||_{L^\infty})d\tau$.
\end{lemma2}

\begin{lemma2}
\cite{BCD} (Existence and uniqueness) Let $p,r,s,f_0$ and $F$ be as in the statement of Lemma 2.1. Assume that
$v\in L^\rho(0,T; B^{-M}_{\infty,\infty})$ for some $\rho >1$ and $M>0$, and
$\partial_x v \in L^1(0,T; {B^{s-1}_{p,r}})$ if $s> 1+{1\over p}$ or $s=1+{1\over p}$ and $r=1$,
and $\partial_x v \in L^1(0,T; B^{1\over p}_{p,\infty}\cap L^\infty)$ if $s<1+{1\over p}$.
Then (T) has a unique solution
$f\in L^{\infty}(0,T; B^s_{p,r})\bigcap\,\big(\bigcap\limits_{s'<s} C([0,T]; B^{s'}_{p,1})\big)$
and the inequalities of Lemma 2.1 can hold true. Moreover, if $r<\infty$, then $f\in C([0,T]; B^s_{p,r})$.
\end{lemma2}

\section{Local well-posedness}
\newtheorem {remark3}{Remark}[section]
\newtheorem{theorem3}{Theorem}[section]
\newtheorem{lemma3}{Lemma}[section]
\newtheorem{corollary3}{Corollary}[section]
\newtheorem{proposition3}{Proposition}[section]
\par
In this section, we will establish the local well-posedness for System (\ref{mn}). To begin with, we apply the classical Kato's semigroup theory \cite{Kato} to obtain the local well-posedness of System (\ref{mn}) in Sobolev spaces.
More precisely, we have
\begin{theorem3}
Suppose that $(m_0, n_0)\in H^{s}(\mathbb{R})\times H^{s}(\mathbb{R})$ with $s\geq 1$.
There exists a maximal existence time $T=T(||m_0||_{H^{s}(\mathbb{R})}, ||n_0||_{H^{s}(\mathbb{R})})>0$,
and a unique solution $(m, n)$ to System (\ref{mn}) such that
$$
(m, n)\in C([0,T);H^{s}(\mathbb{R})\times H^{s}(\mathbb{R}))
\cap C^{1}([0,T);H^{s-1}(\mathbb{R})\times H^{s-1}(\mathbb{R})).
$$
Moreover, the solution depends continuously on the initial data,
that is, the mapping $(m_0, n_0)\mapsto (m, n):$
$$ H^{s}(\mathbb{R})\!\times\! H^{s}(\mathbb{R})\!\rightarrow\!
C([0,T); H^{s}(\mathbb{R})\!\times\! H^{s}(\mathbb{R}))\cap
C^{1}([0,T);H^{s-1}(\mathbb{R})\!\times\! H^{s-1}(\mathbb{R}))
$$
is continuous.
\end{theorem3}

\begin{proof}
By going along the similar line of the proof in \cite{E-Le-Y}, one can readily prove the theorem. For the sake of simplicity, we omit the details here.
\end{proof}

Now we pay attention to the case in the nonhomogeneous Besov spaces. Uniqueness and continuity with respect to the initial data in some sense can be obtained by the following a priori estimates.
\begin{lemma3}
Let $1\leq p,r\leq \infty$ and $s>\max (1-\frac{1}{p},\frac{1}{p},\frac 1 2)$. Suppose that we are given
$(m^{(i)},n^{(i)})\in L^{\infty}(0,T; B^s_{p,r}\times B^s_{p,r})\cap C([0,T];\mathcal{S'}\times \mathcal{S'})$ ($i=1,2$) two solutions of System (\ref{mn}) with the initial data
$(m_0^{(i)},n_0^{(i)})\in B^s_{p,r}\times B^{s}_{p,r}$ ($i=1,2$)
and let $u^{(12)}\triangleq u^{(2)}-u^{(1)}$, $v^{(12)}\triangleq v^{(2)}-v^{(1)}$, and then
$m^{(12)}\triangleq m^{(2)}-m^{(1)}$, $n^{(12)}\triangleq n^{(2)}-n^{(1)}$.
Then for all $t\in[0,T]$, we have\\
(1) if $s>\max (1-\frac{1}{p},\frac{1}{p},\frac 1 2)$, but $s\neq 2+\frac{1}{p}$, then
\begin{eqnarray}
&&||m^{(12)}(t)||_{B^{s-1}_{p,r}}+||n^{(12)}(t)||_{B^{s-1}_{p,r}}\\
\nonumber&\leq&(||m^{(12)}_0||_{B^{s-1}_{p,r}}+||n^{(12)}_0||_{B^{s-1}_{p,r}})\\
\nonumber&&\times e^{C\int_0^t
(||m^{(1)}(\tau)||_{B^s_{p,r}}+||m^{(2)}(\tau)||_{B^s_{p,r}}+||n^{(1)}(\tau)||_{B^{s}_{p,r}}
+||n^{(2)}(\tau)||_{B^{s}_{p,r}}+1)^2 d\tau}\\
\nonumber&\triangleq& L(s-1;t);
\end{eqnarray}
(2) if $s=2+\frac{1}{p}$, then
\begin{eqnarray*}
&& ||m^{(12)}(t)||_{B^{s-1}_{p,r}}+||n^{(12)}(t)||_{B^{s-1}_{p,r}}\\
\nonumber&\leq& C L^{\theta}(s-1;t)((||m^{(1)}(t)||_{B^s_{p,r}}+||m^{(2)}(t)||_{B^s_{p,r}})^{1-\theta}\\
\nonumber&&+(||n^{(1)}(t)||_{B^s_{p,r}}+||n^{(2)}(t)||_{B^s_{p,r}})^{1-\theta}),
\end{eqnarray*}
where $\theta\in(0,1)$.
\end{lemma3}

\begin{proof}
It is obvious that
$(m^{(12)},n^{(12)})\in L^{\infty}(0,T; B^s_{p,r}\times B^s_{p,r})\cap C([0,T];\mathcal{S'}\times \mathcal{S'})$
solves the following Cauchy problem of the transport equations:
\begin{equation}
\left\{\begin{array}{ll}
\partial_t m^{(12)}+\frac{1}{2}(u^{(1)}v^{(1)}-u_x^{(1)}v_x^{(1)})\partial_x m^{(12)}=F(t,x),\\
\partial_t n^{(12)}+\frac{1}{2}(u^{(1)}v^{(1)}-u_x^{(1)}v_x^{(1)})\partial_x n^{(12)}=G(t,x),\\
m^{(12)}|_{t=0}=m^{(12)}_0\triangleq m^{(2)}_0-m^{(1)}_0,\\
n^{(12)}|_{t=0}=n^{(12)}_0\triangleq n^{(2)}_0-n^{(1)}_0,
\end{array}\right.
\end{equation}
where
$F(t,x)\triangleq \frac{1}{2}\big((u^{(1)}-m^{(1)})v_x^{(1)}-(v^{(1)}+n^{(1)})u_x^{(1)}-v_x^{(2)}m^{(2)}\big)m^{(12)}
-\frac{1}{2}u_x^{(2)}m^{(2)}n^{(12)}
+\frac{1}{2}\big(v_x^{(1)}m_x^{(2)}-(v^{(1)}+n^{(1)})m^{(2)}-2b\big)u_x^{(12)}
+\frac{1}{2}\big(u_x^{(2)}m_x^{(2)}-(m^{(1)}-u^{(2)})m^{(2)}\big)v_x^{(12)}
+\frac{1}{2}\big(v_x^{(1)}m^{(2)}-v^{(1)}m_x^{(2)}\big)u^{(12)}
-\frac{1}{2}(u^{(2)}m^{(2)})_x v^{(12)}$, \,\,and\\
$G(t,x)\triangleq \frac{1}{2}\big((v^{(1)}-n^{(1)})u_x^{(1)}-(u^{(1)}+m^{(1)})v_x^{(1)}-u_x^{(2)}n^{(2)}\big)n^{(12)}
-\frac{1}{2}v_x^{(2)}n^{(2)}m^{(12)}
+\frac{1}{2}\big(u_x^{(1)}n_x^{(2)}-(u^{(1)}+m^{(1)})n^{(2)}-2b\big)v_x^{(12)}
+\frac{1}{2}\big(v_x^{(2)}n_x^{(2)}-(n^{(1)}-v^{(2)})n^{(2)}\big)u_x^{(12)}
+\frac{1}{2}\big(u_x^{(1)}n^{(2)}-u^{(1)}n_x^{(2)}\big)v^{(12)}
-\frac{1}{2}(v^{(2)}n^{(2)})_x u^{(12)}$.\\

\noindent {\it{Claim}}. For all $s>\max (\frac{1}{p},\frac 1 2)$ and $t\in [0,T]$, we have
\begin{eqnarray}
&&||F(t)||_{B^{s-1}_{p,r}},\,\ ||G(t)||_{B^{s-1}_{p,r}}\\
\nonumber&\leq&C(||m^{(12)}(t)||_{B^{s-1}_{p,r}}+||n^{(12)}(t)||_{B^{s-1}_{p,r}})\\
\nonumber&&\times (||m^{(1)}(t)||_{B^s_{p,r}}+||m^{(2)}(t)||_{B^s_{p,r}}+||n^{(1)}(t)||_{B^{s}_{p,r}}
+||n^{(2)}(t)||_{B^{s}_{p,r}}+1)^2,
\end{eqnarray}
where $C=C(s,p,r,b)$ is a positive constant.\\
Indeed, for $s>1+\frac{1}{p}$, $B^{s-1}_{p,r}$ is an algebra, by Proposition 2.2 (4), we have
\begin{eqnarray*}
||(u^{(1)}-m^{(1)})v_x^{(1)}m^{(12)}||_{B^{s-1}_{p,r}}
\leq(||u^{(1)}||_{B^{s-1}_{p,r}}+||m^{(1)}||_{B^{s-1}_{p,r}})
||v_x^{(1)}||_{B^{s-1}_{p,r}}||m^{(12)}||_{B^{s-1}_{p,r}}.
\end{eqnarray*}
Note that $(1-\partial_x^2)^{-1}\in Op(S^{-2})$. According to Proposition 2.2 (8), we obtain
for all $s\in\mathbb{R}$ and $i=1,2,12$,
\begin{eqnarray}
||u^{(i)}||_{B^{s+2}_{p,r}}\approxeq ||m^{(i)}||_{B^{s}_{p,r}}\quad  \text{and}\quad
||v^{(i)}||_{B^{s+2}_{p,r}}\approxeq ||n^{(i)}||_{B^{s}_{p,r}}.
\end{eqnarray}
Then
\begin{eqnarray*}
||(u^{(1)}-m^{(1)})v_x^{(1)}m^{(12)}||_{B^{s-1}_{p,r}}
\leq C||m^{(1)}||_{B^{s}_{p,r}}||n^{(1)}||_{B^{s}_{p,r}} ||m^{(12)}||_{B^{s-1}_{p,r}}.
\end{eqnarray*}
Similarly,
\begin{eqnarray*}
&&||(v^{(1)}+n^{(1)})u_x^{(1)}m^{(12)}||_{B^{s-1}_{p,r}}+||v_x^{(2)}m^{(2)}m^{(12)}||_{B^{s-1}_{p,r}}\\ \nonumber
&\leq& C (||m^{(1)}||_{B^{s}_{p,r}}||n^{(1)}||_{B^{s}_{p,r}}+ ||m^{(2)}||_{B^{s}_{p,r}}||n^{(2)}||_{B^{s}_{p,r}})||m^{(12)}||_{B^{s-1}_{p,r}},
\end{eqnarray*}
\begin{eqnarray*}
||u_x^{(2)}m^{(2)}n^{(12)}||_{B^{s-1}_{p,r}}+||(u^{(2)}m^{(2)})_x v^{(12)}||_{B^{s-1}_{p,r}}
\leq C ||m^{(2)}||^2_{B^{s}_{p,r}}||n^{(12)}||_{B^{s-1}_{p,r}},
\end{eqnarray*}
\begin{eqnarray*}
&&||\big(v_x^{(1)}m_x^{(2)}-(v^{(1)}+n^{(1)})m^{(2)}-2b\big)u_x^{(12)}||_{B^{s-1}_{p,r}}
+||\big(v_x^{(1)}m^{(2)}-v^{(1)}m_x^{(2)}\big)u^{(12)}||_{B^{s-1}_{p,r}}\\ \nonumber
&\leq& C (||m^{(2)}||_{B^{s}_{p,r}}||n^{(1)}||_{B^{s}_{p,r}}+ 1)||m^{(12)}||_{B^{s-1}_{p,r}},
\end{eqnarray*}
and
\begin{eqnarray*}
&&||\big(u_x^{(2)}m_x^{(2)}-(m^{(1)}-u^{(2)})m^{(2)}\big)v_x^{(12)}||_{B^{s-1}_{p,r}}\\ \nonumber
&\leq& C (||m^{(1)}||_{B^{s}_{p,r}}+||m^{(2)}||_{B^{s}_{p,r}})||m^{(2)}||_{B^{s}_{p,r}}||n^{(12)}||_{B^{s-1}_{p,r}}.
\end{eqnarray*}
Hence,
\begin{eqnarray*}
||F(t)||_{B^{s-1}_{p,r}} &\leq& C(||m^{(12)}(t)||_{B^{s-1}_{p,r}}+||n^{(12)}(t)||_{B^{s-1}_{p,r}})\\ \nonumber
&\times& (||m^{(1)}(t)||_{B^s_{p,r}}+||m^{(2)}(t)||_{B^s_{p,r}}+||n^{(1)}(t)||_{B^{s}_{p,r}}
+||n^{(2)}(t)||_{B^{s}_{p,r}}+1)^2,
\end{eqnarray*}
if $s>1+\frac{1}{p}$.\\
We can treat $||G(t)||_{B^{s-1}_{p,r}}$ for $s>1+\frac{1}{p}$ in a similar way.

On the other hand, if $\max (\frac{1}{p},\frac 1 2)<s\leq {1+\frac{1}{p}}$, $B^{s}_{p,r}$ is an algebra. According to (2.1) and (3.4), one infers that
\begin{eqnarray*}
&&||\big((u^{(1)}-m^{(1)})v_x^{(1)}-(v^{(1)}+n^{(1)})u_x^{(1)}-v_x^{(2)}m^{(2)}\big)m^{(12)}||_{B^{s-1}_{p,r}}\\
&\leq& C (||(u^{(1)}-m^{(1)})v_x^{(1)}||_{B^{s}_{p,r}}+||(v^{(1)}+n^{(1)})u_x^{(1)}||_{B^{s}_{p,r}}+
||v_x^{(2)}m^{(2)}||_{B^{s}_{p,r}})||m^{(12)}||_{B^{s-1}_{p,r}}\\
&\leq& C (||m^{(1)}||_{B^s_{p,r}}||n^{(1)}||_{B^s_{p,r}}+||m^{(2)}||_{B^s_{p,r}}||n^{(2)}||_{B^s_{p,r}})
||m^{(12)}||_{B^{s-1}_{p,r}},
\end{eqnarray*}
\begin{eqnarray*}
||u_x^{(2)}m^{(2)}n^{(12)}||_{B^{s-1}_{p,r}}
\leq C ||u_x^{(2)}m^{(2)}||_{B^{s}_{p,r}}||n^{(12)}||_{B^{s-1}_{p,r}}
\leq C ||m^{(2)}||^2_{B^{s}_{p,r}}||n^{(12)}||_{B^{s-1}_{p,r}},
\end{eqnarray*}
\begin{eqnarray*}
&&||\big(v_x^{(1)}m_x^{(2)}-(v^{(1)}+n^{(1)})m^{(2)}-2b\big)u_x^{(12)}||_{B^{s-1}_{p,r}}\\
&\leq& C ||v_x^{(1)}u_x^{(12)}||_{B^{s}_{p,r}}||m_x^{(2)}||_{B^{s-1}_{p,r}}
+C (||(v^{(1)}+n^{(1)})m^{(2)}||_{B^{s}_{p,r}}+|b|)||u_x^{(12)}||_{B^{s-1}_{p,r}}\\
&\leq& C (||m^{(2)}||_{B^s_{p,r}}||n^{(1)}||_{B^s_{p,r}}+1)||u^{(12)}||_{B^{s+1}_{p,r}}\\
&\leq& C (||m^{(2)}||_{B^s_{p,r}}||n^{(1)}||_{B^s_{p,r}}+1)||m^{(12)}||_{B^{s-1}_{p,r}},
\end{eqnarray*}
\begin{eqnarray*}
&&||\big(u_x^{(2)}m_x^{(2)}-(m^{(1)}-u^{(2)})m^{(2)}\big)v_x^{(12)}||_{B^{s-1}_{p,r}}\\
&\leq& C (||u_x^{(2)}v_x^{(12)}||_{B^{s}_{p,r}}||m_x^{(2)}||_{B^{s-1}_{p,r}}
+||(m^{(1)}-u^{(2)})m^{(2)}||_{B^{s}_{p,r}}||v_x^{(12)}||_{B^{s-1}_{p,r}})\\
&\leq& C (||m^{(1)}||_{B^s_{p,r}}+||m^{(2)}||_{B^s_{p,r}})||m^{(2)}||_{B^s_{p,r}}||v^{(12)}||_{B^{s+1}_{p,r}}\\
&\leq& C (||m^{(1)}||_{B^s_{p,r}}+||m^{(2)}||_{B^s_{p,r}})||m^{(2)}||_{B^s_{p,r}}||n^{(12)}||_{B^{s-1}_{p,r}},
\end{eqnarray*}
\begin{eqnarray*}
&&||\big(v_x^{(1)}m^{(2)}-v^{(1)}m_x^{(2)}\big)u^{(12)}||_{B^{s-1}_{p,r}}\\
&\leq& C (||v_x^{(1)}m^{(2)}||_{B^{s}_{p,r}}||u^{(12)}||_{B^{s-1}_{p,r}}
+||v^{(1)}u^{(12)}||_{B^{s}_{p,r}}||m_x^{(2)}||_{B^{s-1}_{p,r}})\\
&\leq& C ||m^{(2)}||_{B^s_{p,r}}||n^{(1)}||_{B^s_{p,r}}||m^{(12)}||_{B^{s-1}_{p,r}},
\end{eqnarray*}
and
\begin{eqnarray*}
&&||(u^{(2)}m^{(2)})_x v^{(12)}||_{B^{s-1}_{p,r}}\\
&\leq& C (||u_x^{(2)}m^{(2)}||_{B^{s}_{p,r}}||v^{(12)}||_{B^{s-1}_{p,r}}
+||u^{(2)}v^{(12)}||_{B^{s}_{p,r}}||m_x^{(2)}||_{B^{s-1}_{p,r}})\\
&\leq& C ||m^{(2)}||^2_{B^s_{p,r}}||n^{(12)}||_{B^{s-1}_{p,r}}.
\end{eqnarray*}
Thus,
\begin{eqnarray*}
||F(t)||_{B^{s-1}_{p,r}} &\leq& C(||m^{(12)}(t)||_{B^{s-1}_{p,r}}+||n^{(12)}(t)||_{B^{s-1}_{p,r}})\\ \nonumber
&\times& (||m^{(1)}(t)||_{B^s_{p,r}}+||m^{(2)}(t)||_{B^s_{p,r}}+||n^{(1)}(t)||_{B^{s}_{p,r}}
+||n^{(2)}(t)||_{B^{s}_{p,r}}+1)^2,
\end{eqnarray*}
if $\max (\frac{1}{p},\frac 1 2)<s\leq {1+\frac{1}{p}}$.\\
We can also treat $||G(t)||_{B^{s-1}_{p,r}}$ for $\max (\frac{1}{p},\frac 1 2)<s\leq {1+\frac{1}{p}}$ in a similar way. Therefore, we prove our Claim (3.3).

Applying Lemma 2.1 (1) and the fact that
\begin{eqnarray*}
V(t)&\triangleq&||\partial_x(u^{(1)}v^{(1)}-u_x^{(1)}v_x^{(1)})||_{B^{\frac 1 p}_{p,r}\cap L^{\infty}}
+||\partial_x(u^{(1)}v^{(1)}-u_x^{(1)}v_x^{(1)})||_{B^{s-2}_{p,r}}\\
&\leq& C ||\partial_x(u^{(1)}v^{(1)}-u_x^{(1)}v_x^{(1)})||_{B^{s}_{p,r}}\\
&\leq& C ||u^{(1)}||_{B^{s+2}_{p,r}}||v^{(1)}||_{B^{s+2}_{p,r}}\\
&\leq& C ||m^{(1)}||_{B^s_{p,r}}||n^{(1)}||_{B^{s}_{p,r}}
\end{eqnarray*}
for $s>\max (1-\frac{1}{p},\frac{1}{p},\frac 1 2)$, we can obtain, for case (1),
$$||m^{(12)}(t)||_{B^{s-1}_{p,r}}\leq ||m^{(12)}_0||_{B^{s-1}_{p,r}}+\int_0^t ||F(\tau)||_{B^{s-1}_{p,r}}d\tau +C \int_0^t V(\tau)||m^{(12)}(\tau)||_{B^{s-1}_{p,r}}d\tau$$
and
$$||n^{(12)}(t)||_{B^{s-1}_{p,r}}\leq ||n^{(12)}_0||_{B^{s-1}_{p,r}}+\int_0^t ||G(\tau)||_{B^{s-1}_{p,r}}d\tau +C \int_0^t V(\tau)||n^{(12)}(\tau)||_{B^{s-1}_{p,r}}d\tau$$
which together with (3.3) yields
\begin{eqnarray*}
&&||m^{(12)}(t)||_{B^{s-1}_{p,r}}+||n^{(12)}(t)||_{B^{s-1}_{p,r}}\\
\nonumber&\leq&||m^{(12)}_0||_{B^{s-1}_{p,r}}+||n^{(12)}_0||_{B^{s-1}_{p,r}}+C \int_0^t (||m^{(12)}(\tau)||_{B^{s-1}_{p,r}}+||n^{(12)}(\tau)||_{^{s-1}_{p,r}})\\
\nonumber&&\times (||m^{(1)}(\tau)||_{B^s_{p,r}}+||m^{(2)}(\tau)||_{B^s_{p,r}}
+||n^{(1)}(\tau)||_{B^s_{p,r}}+||n^{(2)}(\tau)||_{B^s_{p,r}}+1)^2 d\tau.
\end{eqnarray*}
Taking advantage of Gronwall's inequality, we get (3.1).

For the critical case (2) $s= 2+\frac{1}{p}$, we here use the interpolation method to deal with it. Indeed, if we choose
$s_1\in(\max(1-\frac{1}{p},\frac{1}{p},\frac 1 2)-1,s-1)$, $s_2\in (s-1,s)$ and $\theta=\frac{s_2-(s-1)}{s_2-s_1} \in(0,1)$,
then $s-1=\theta s_1+(1-\theta) s_2$. According to Proposition 2.2 (6) and the consequence of case (1), we have
\begin{eqnarray*}
&&||m^{(12)}(t)||_{B^{s-1}_{p,r}}+||n^{(12)}(t)||_{B^{s-1}_{p,r}}\\
\nonumber&\leq&||m^{(12)}(t)||^{\theta}_{B^{s_1}_{p,r}}||m^{(12)}(t)||^{1-\theta}_{B^{s_2}_{p,r}}
+||n^{(12)}(t)(t)||^{\theta}_{B^{s_1}_{p,r}}||n^{(12)}(t)(t)||^{1-\theta}_{B^{s_2}_{p,r}}\\
\nonumber&\leq& C L^{\theta}(s_1;t)((||m^{(1)}(t)||_{B^{s_2}_{p,r}}+||m^{(2)}(t)||_{B^{s_2}_{p,r}})^{1-\theta}\\
\nonumber&&+(||n^{(1)}(t)||_{B^{s_2}_{p,r}}+||n^{(2)}(t)||_{B^{s_2}_{p,r}})^{1-\theta})\\
\nonumber&\leq& C L^{\theta}(s-1;t)((||m^{(1)}(t)||_{B^s_{p,r}}+||m^{(2)}(t)||_{B^s_{p,r}})^{1-\theta}\\
\nonumber&&+(||n^{(1)}(t)||_{B^s_{p,r}}+||n^{(2)}(t)||_{B^s_{p,r}})^{1-\theta})
\end{eqnarray*}
Therefore, we complete our proof of Lemma 3.1.
\end{proof}

We next construct the smooth approximation solutions to System (\ref{mn}) as follows.
\begin{lemma3}
Let $p$ and $r$ be as in the statement of Lemma 3.1. Assume that $s>\max (1-\frac{1}{p},\frac{1}{p},\frac 1 2)$ and
$s\neq1+\frac{1}{p}$, $(m_0,n_0)\in B^s_{p,r}\times B^{s}_{p,r}$ and $(m^0,n^0)=(0,0)$. Then
(1) there exists a sequence of smooth functions $(m^k,n^k)_{k\in \mathbb{N}}$ belonging to
$C(\mathbb{R}^{+}; B^{\infty}_{p,r}\times B^{\infty}_{p,r})$
and solving the following linear transport equations by induction with respect to $k$:
\[(T_k)\left\{\begin{array}{l}
\partial_t m^{k+1}+\frac{1}{2}(u^k v^k-u^k_x v^k_x)\partial_x m^{k+1}=R^k_1(t,x),\\
\partial_t n^{k+1}+\frac{1}{2}(u^k v^k-u^k_x v^k_x)\partial_x n^{k+1}=R^k_2(t,x),\\
m^{k+1}|_{t=0}\triangleq m^{k+1}_0(x)=S_{k+1}m_0,\\
n^{k+1}|_{t=0}\triangleq n^{k+1}_0(x)=S_{k+1}n_0,
\end{array}\right.\]
where $R^k_1(t,x)\triangleq -\frac{1}{2}\big((u^k v^k-u^k_x v^k_x)_x-(u^k v^k_x-u^k_x v^k)\big)m^k-bu^k_x$ and
$R^k_2(t,x)\triangleq -\frac{1}{2}\big((u^k v^k-u^k_x v^k_x)_x+(u^k v^k_x-u^k_x v^k)\big)n^k-bv^k_x$.\\
(2) there exists $T>0$ such that the solution $(m^k,n^k)_{k\in \mathbb{N}}$
is uniformly bounded in $E^s_{p,r}(T)\times E^{s}_{p,r}(T)$ and a Cauchy sequence in
$C([0,T]; B^{s-1}_{p,r}\times B^{s-1}_{p,r})$, whence it converges to some limit
$(m,n)\in C([0,T]; B^{s-1}_{p,r}\times B^{s-1}_{p,r})$.
\end{lemma3}

\begin{proof}
Since all the data $S_{k+1}m_0,\,S_{k+1}n_0\in B^{\infty}_{p,r}$, it then follows from Lemma 2.3 and by induction with respect to the index $k$ that (1) holds.

To prove (2), applying Remark (2.2) and  simulating  the proof of Lemma 3.1 (1), we obtain that for
$s>\max (1-\frac{1}{p},\frac{1}{p},\frac 1 2)$ and $s\neq1+\frac{1}{p}$,
\begin{eqnarray}
&& a_{k+1}(t)\\ \nonumber
&\leq& C e^{C U^k(t)}\big(A+\int_0^t e^{-C U^k(\tau)}
(||R^k_1(\tau)||_{B^{s}_{p,r}}+||R^k_2(\tau)||_{B^{s}_{p,r}})d\tau\big),
\end{eqnarray}
where $a_k(t)\triangleq ||m^k(t)||_{B^s_{p,r}}+||n^k(t)||_{B^{s}_{p,r}}$,
$A\triangleq ||m_0||_{B^s_{p,r}}+||n_0||_{B^{s}_{p,r}}$
and $U^k(t)\triangleq \int_0^t ||m^k(\tau)||_{B^s_{p,r}}||n^k(\tau)||_{B^s_{p,r}}d\tau$.\\
Noting that $B^{s}_{p,r}$ is an algebra and according to (3.4), one gets
\begin{eqnarray*}
&&||R^k_1(t)||_{B^{s}_{p,r}}+||R^k_2(t)||_{B^{s}_{p,r}}\\ \nonumber
&\leq& C (||m^k(t)||_{B^s_{p,r}}+||n^k(t)||_{B^{s}_{p,r}})
(1+||m^k(t)||_{B^s_{p,r}}||n^k(t)||_{B^{s}_{p,r}})\\ \nonumber
&\leq& C (a_k(t)+a^3_k(t)).
\end{eqnarray*}
If $a_k(t)<1$, then by (3.5), $a_{k+1}(t)\leq C(A+t)$, which implies that $(m^k,n^k)_{k\in \mathbb{N}}$
is uniformly bounded in $C([0,T]; B^{s}_{p,r}\times B^{s}_{p,r})$.\\
If $a_k(t)\geq1$, from (3.5), we have
\begin{eqnarray}
a_{k+1}(t) \leq C e^{C U^k(t)}\big(A+\int_0^t e^{-C U^k(\tau)}a^3_k(\tau)d\tau\big).
\end{eqnarray}
Choose $0<\, T\, < \frac{1}{4C^3 A^2}$ and suppose that
\begin{eqnarray}
a_k(t)\leq \frac{CA}{\sqrt{1-4C^3 A^2t}},\quad \forall t\in[0,T].
\end{eqnarray}
Noting that $e^{C (U^k(t)-U^k(\tau))}\leq \sqrt[4]{\frac{1-4C^3 A^2\tau}{1-4C^3 A^2t}}$ and substituting (3.7) into (3.6) yields
\begin{eqnarray*}
a_{k+1}(t)
&\leq&\frac{CA}{\sqrt[4]{1-4C^3 A^2t}}+\frac{C}{\sqrt[4]{1-4C^3 A^2t}}
\int_0^t \frac{C^3 A^3}{(1-4C^3 A^2\tau)^\frac{5}{4}} d\tau\\ \nonumber
&=&\frac{CA}{\sqrt[4]{1-4C^3 A^2t}}+\frac{C}{\sqrt[4]{1-4C^3 A^2t}} (\frac{A}{\sqrt[4]{1-4C^3 A^2t}}-A)\\ \nonumber
&\leq&\frac{CA}{\sqrt{1-4C^3 A^2t}},
\end{eqnarray*}
which implies that
$$(m^k,n^k)_{k\in \mathbb{N}}\quad \text{is uniformly bounded in}\quad C([0,T]; B^{s}_{p,r}\times B^{s}_{p,r}).$$
Using the equations $(T_k)$ and the similar argument in the proof of Lemma 3.1 (1), one can easily prove that
$$(\partial_t m^{k+1},\partial_t m^{k+1})_{k\in \mathbb{N}}\quad
\text{is uniformly  bounded in} \quad  C([0,T]; B^{s-1}_{p,r}\times B^{s-1}_{p,r}).$$
Hence,
$$(m^k,n^k)_{k\in \mathbb{N}}\quad \text{is uniformly bounded in}\quad E^s_{p,r}(T)\times E^{s}_{p,r}(T).$$

Now it suffices to show that $(m^k,n^k)_{k\in \mathbb{N}}$
is a Cauchy sequence in $C([0,T]; B^{s-1}_{p,r})\times C([0,T]; B^{s-1}_{p,r})$.
Indeed, for all $k,l\in \mathbb{N}$, from $(T_k)$, we have
\begin{eqnarray*}
&&\partial_t (m^{k+l+1}-m^{k+1})+\frac{1}{2}(u^{k+l}v^{k+l}-u_x^{k+l}v_x^{k+l}) \partial_x (m^{k+l+1}-m^{k+1})\\
\nonumber&=&
-\frac{1}{2}\big((u^{k+l}v^{k+l}-u_x^{k+l}v_x^{k+l})_x-(u^{k+l}v_x^{k+l}-u_x^{k+l}v^{k+l})\big)(m^{k+l}-m^k)\\ \nonumber&&-\frac{1}{2}
\big(
(u^{k+l}v^{k+l}-u^k v^k)_x-(u_x^{k+l}v_x^{k+l}-u_x^k v_x^k)_x-(u^{k+l}v_x^{k+l}-u^k v_x^k)\\ \nonumber
\nonumber&&+(u_x^{k+l}v^{k+l}-u_x^k v^k)
\big)m^k
-\frac{1}{2}\big((u^{k+l}v^{k+l}-u^k v^k)-(u_x^{k+l}v_x^{k+l}-u_x^k v_x^k) \big)m_x^{k+1}\\ \nonumber
\nonumber&& -b (u_x^{k+l}-u_x^k)
\end{eqnarray*}
and
\begin{eqnarray*}
&&\partial_t (n^{k+l+1}-n^{k+1})+\frac{1}{2}(u^{k+l}v^{k+l}-u_x^{k+l}v_x^{k+l}) \partial_x (n^{k+l+1}-n^{k+1})\\
\nonumber&=&
-\frac{1}{2}\big((u^{k+l}v^{k+l}-u_x^{k+l}v_x^{k+l})_x+(u^{k+l}v_x^{k+l}-u_x^{k+l}v^{k+l})\big)(n^{k+l}-n^k)\\ \nonumber&&-\frac{1}{2}
\big(
(u^{k+l}v^{k+l}-u^k v^k)_x-(u_x^{k+l}v_x^{k+l}-u_x^k v_x^k)_x-(u_x^{k+l}v^{k+l}-u_x^k v^k)\\ \nonumber
\nonumber&&+(u^{k+l}v_x^{k+l}-u^k v_x^k)
\big)n^k
-\frac{1}{2}\big((u^{k+l}v^{k+l}-u^k v^k)-(u_x^{k+l}v_x^{k+l}-u_x^k v_x^k) \big)n_x^{k+1}\\ \nonumber
\nonumber&& -b (v_x^{k+l}-v_x^k).
\end{eqnarray*}
Similar to the proof of Lemma 3.1 (1), for $s>\max (1-\frac{1}{p},\frac{1}{p},\frac 1 2)$ and $s\neq 1+\frac{1}{p}, 2+\frac{1}{p}$, we can obtain that
$$b^l_{k+1}(t)\leq C e^{C U^{k+l}(t)}\big(b^l_{k+1}(0)+ \int_0^t
e^{-C U^{k+l}(\tau)}d^l_k(\tau)b^l_k(\tau)d\tau\big),$$
where $b^l_k(t)\triangleq ||(m^{k+l}-m^k)(t)||_{B^{s-1}_{p,r}}+||(n^{k+l}-n^k)(t)||_{B^{s-1}_{p,r}}$,\\ $U^{k+l}(t)\triangleq \int_0^t ||m^{k+l}(\tau)||_{B^s_{p,r}}||m^{k+l}(\tau)||_{B^s_{p,r}} d\tau$, \quad
and\\
$d^l_k(t)\triangleq (||m^k(t)||_{B^s_{p,r}}+||m^{k+1}(t)||_{B^s_{p,r}}+||m^{k+l}(t)||_{B^s_{p,r}}
+||n^k(t)||_{B^s_{p,r}}+||n^{k+1}(t)||_{B^s_{p,r}}\\
+||n^{k+l}(t)||_{B^s_{p,r}})^2+1$.\\
Thanks to Remark 2.1, we have
\begin{eqnarray*}
||\sum\limits_{q=k+1}^{k+l} \Delta_q m_0||_{B^{s-1}_{p,r}}
&=&
\big(\sum\limits_{j\geq-1}2^{j(s-1)r}||\Delta_j (\sum\limits_{q=k+1}^{k+l} \Delta_q m_0)||^r_{L^p}\big)^{\frac 1 r}\\
&\leq& C \big(\sum\limits_{j=k}^{k+l+1}2^{-jr}2^{jsr}||\Delta_j{m_0}||^r_{L^p}\big)^{\frac 1 r}\\
&\leq& C2^{-k}||m_0||_{B^s_{p,r}}.
\end{eqnarray*}
Similarly,
$$||\sum\limits_{q=k+1}^{k+l} \Delta_q n_0||_{B^{s-1}_{p,r}}\leq C2^{-k}||n_0||_{B^{s}_{p,r}}.$$
Hence, we obtain
$$b^l_{k+1}(0)\leq C2^{-k}(||m_0||_{B^s_{p,r}}+||n_0||_{B^{s}_{p,r}}).$$
According to the fact that $(m^k,n^k)_{k\in \mathbb{N}}$
is uniformly bounded in $E^s_{p,r}(T)\times E^{s}_{p,r}(T)$, we can find a positive constant $C_{T}$ independent of $k,l$ such that
$$b^l_{k+1}(t)\leq C_{T}\big(2^{-k}+\int_0^t b^l_k(\tau)d\tau\big),\quad \forall t\in[0,T].$$
Arguing by induction with respect to the index $k$, we can obtain
\begin{eqnarray*}
b^l_{k+1}(t)
&\leq& C_{T}\big(2^{-k}\sum\limits_{j=0}^k \frac{(2T C_T)^j}{j!}+C^{k+1}_T \int_0^t  \frac{(t-\tau)^k}{k!}d\tau\big)\\
&\leq& \big(C_{T}\sum\limits_{j=0}^k \frac{(2T C_T)^j}{j!}\big)2^{-k}+C_T \frac{(T C_T)^{k+1}}{(k+1)!},
\end{eqnarray*}
which implies the desired result as $k\rightarrow +\infty$.

On the other hand, for the critical point $2+\frac{1}{p}$, we can apply the interpolation method which has been used in the proof of Lemma 3.1 to show that $(m^k,n^k)_{k\in \mathbb{N}}$
is also a Cauchy sequence in $C([0,T]; B^{s-1}_{p,r}\times B^{s-1}_{p,r})$ for this critical case.
Therefore, we have completed the proof of Lemma 3.2.
\end{proof}

Now we are in the position to prove the main theorem of this section.
\begin{theorem3}
Assume that $1\leq p,r\leq \infty$ and $s>\max (1-\frac{1}{p},\frac{1}{p},\frac 1 2)$ but $s\neq1+\frac{1}{p}$. Let
$(m_0,n_0)\in B^s_{p,r}\times B^{s}_{p,r}$ and
$(m,n)$ be the obtained limit in Lemma 3.2 (2).
Then there exists a time $T>0$ such that $(m,n)\in E^s_{p,r}(T)\times E^{s}_{p,r}(T)$ is the unique solution to System (\ref{mn}), and the mapping $(m_0,n_0)\mapsto (m,n)$ is continuous from
$ B^s_{p,r}\times B^{s}_{p,r}$ into
$$ C([0,T]; B^{s'}_{p,r}\times B^{s'}_{p,r})\cap
C^1{([0,T]; B^{s'-1}_{p,r}\times B^{s'-1}_{p,r})}$$
for all $s'<s$ if $r=\infty$ and $s'=s$ otherwise.
\end{theorem3}

\begin{proof}
We first claim that $(m,n)\in E^s_{p,r}(T)\times E^{s}_{p,r}(T)$ solves System (\ref{mn}).\\
In fact, according to Lemma 3.2 (2) and Proposition 2.2 (7), one can get
$$(m,n)\in L^{\infty}([0,T]; B^s_{p,r}\times B^s_{p,r}).$$
For all $s'<s$, Lemma 3.2 (2) applied again, together with an interpolation argument yields
$$(m^k,n^k) \rightarrow (m,n),\,\ \text{as} \,\ n\to \infty,\,\ \text{in}\,\
C([0,T]; B^{s'}_{p,r}\times B^{s'}_{p,r}). $$
Taking limit in $(T_k)$, we can see that
$(m,n)$ solves System (\ref{mn}) in the sense of $C([0,T]; B^{s'-1}_{p,r}\times B^{s'-1}_{p,r})$ for all $s'<s$.\\
Making use of the equations in System (\ref{mn}) twice and the similar proof of (3.4), together with Lemma 2.1 (3) and Lemma 2.3 yields
$(m,n)\in E^s_{p,r}(T)\times E^{s}_{p,r}(T).$

On the other hand, the continuity with respect to the initial data in
$$ C([0,T]; B^{s'}_{p,r}\times B^{s'}_{p,r})\cap C^1{([0,T]; B^{s'-1}_{p,r}\times B^{s'-1}_{p,r})}
\quad (\forall \, s'<s)$$
can be obtained by Lemma 3.1 and a simple interpolation argument. While the continuity in
$C([0,T]; B^{s}_{p,r}\times B^{s}_{p,r})\cap C^1{([0,T]; B^{s-1}_{p,r}\times B^{s-1}_{p,r})}$ when $r< \infty$ can be proved through the use of a sequence of viscosity approximation solutions
$(m_\varepsilon,n_\varepsilon)_{\varepsilon>0}$
for System (\ref{mn}) which converges uniformly in
$C([0,T]; B^{s}_{p,r}\times B^{s}_{p,r})\cap C^1{([0,T]; B^{s-1}_{p,r}\times B^{s-1}_{p,r})}$.
This completes the proof of Theorem 3.2.
\end{proof}

\begin{remark3}
Note that for every $s\in\mathbb{R}$, $B^{s}_{2,2}=H^s$. Theorem 3.2 holds true in the corresponding Sobolev spaces with $\frac{1} {2}<s\neq\frac{3} {2}$, which almost improves the result of Theorem 3.1 proved by Kato's theory, where $s\geq 1$ is required. Therefore, Theorem 3.2 together with Theorem 3.1 implies that the conclusion of Theorem 3.1  holds true for initial data $(m_0,n_0)\in H^s(\mathbb{R})\times H^s(\mathbb{R})$ with $s>\frac{1} {2}$
or for all initial data $(u_0,v_0)\in H^s(\mathbb{R})\times H^s(\mathbb{R})$ with $s>\frac{5} {2}$. Moreover, by going along the similar proof of Theorems 3.1-3.2, one can readily obtain that Theorem 3.1 is also true for System (\ref{mn}) in terms of the solution $(u+C_1,v+C_2)$ related to the initial data
$(u_0+C_1,v_0+C_2)\in H^s(\mathbb{R})\times H^s(\mathbb{R})$ with $s>\frac{5} {2}$ for two constants $C_1$ and $C_2$.
\end{remark3}

\section{Blow-up}
\newtheorem {remark4}{Remark}[section]
\newtheorem{theorem4}{Theorem}[section]
\newtheorem{lemma4}{Lemma}[section]
\newtheorem{corollary4}{Corollary}[section]

In this section, we will derive the precise blow-up scenario of strong solutions to System (\ref{mn})
and then state a new blow-up result with respect to initial data.

\begin{theorem4}
Let $(m_0,n_0)\in H^s(\mathbb{R})\times H^{s}(\mathbb{R})$ with $s>\frac 1 2$ and
$T$ be the maximal existence time of the solution $(m,n)$ to System (\ref{mn}), which is guaranteed by Remark 3.1.
If $T<\infty$, then
$$\int_0^T ||m(\tau)||_{L^\infty} ||n(\tau)||_{L^\infty} d\tau=\infty.$$
\end{theorem4}

\begin{proof}
 We will prove the theorem by induction with respect to the regular index $s$ $(s>\frac 1 2)$ as follows.

\noindent {\it{Step 1}}. For $s\in(\frac{1}{2},1)$, by Lemma 2.2 and System (\ref{mn}), we have
\begin{eqnarray*}
||m(t)||_{H^{s}}&\leq& ||m_0||_{H^{s}}
+C\int_0^t ||m(\tau)||_{H^s}(||uv-u_x v_x||_{L^{\infty}}+||u_x n+v_x m||_{L^{\infty}}) d\tau\\
&&+C \int_0^t
(||\big((u_x n+v_x m)-(uv_x-u_xv)\big)m(\tau)||_{H^s}+||u_x(\tau)||_{H^s})d\tau
\end{eqnarray*}
and
\begin{eqnarray*}
||n(t)||_{H^{s}}&\leq& ||n_0||_{H^{s}}
+C\int_0^t ||n(\tau)||_{H^s}(||uv-u_x v_x||_{L^{\infty}}+||u_x n+v_x m||_{L^{\infty}}) d\tau\\
&&+C \int_0^t
(||\big((u_x n+v_x m)+(uv_x-u_x v)\big)n(\tau)||_{H^s}+||v_x(\tau)||_{H^s})d\tau.
\end{eqnarray*}
Noting that $u=(1-\partial_x^2)^{-1}m=p\ast m$ with $p(x)\triangleq \frac{1}{2}e^{-|x|}\,(x\in\mathbb{R})$, $u_x=\partial_x p\ast m$, $u_{xx}=u-m$ and $||p||_{L^1}=||\partial_x p||_{L^1}=1$, together with the Young inequality implies that for all $s\in \mathbb{R}$,
\begin{eqnarray}
||u||_{L^\infty},\,||u_x||_{L^\infty},\,||u_{xx}||_{L^\infty}\leq C ||m||_{L^\infty}
\end{eqnarray}
and
\begin{eqnarray}
||u||_{H^s},\,||u_x||_{H^s},\,||u_{xx}||_{H^s}\leq C ||m||_{H^s}.
\end{eqnarray}
Similarly, the identity $v=p\ast n$ ensures
\begin{eqnarray}
||v||_{L^\infty},\,||v_x||_{L^\infty},\,||v_{xx}||_{L^\infty}\leq C ||n||_{L^\infty}
\end{eqnarray}
and
\begin{eqnarray}
||v||_{H^s},\,||v_x||_{H^s},\,||v_{xx}||_{H^s}\leq C ||n||_{H^s}.
\end{eqnarray}
Then
\begin{eqnarray}
&&||\big((u_x n+v_x m)-(uv_x-u_xv)\big)m||_{H^s}+||u_x||_{H^s}\\ \nonumber
&\leq& C (||m||_{L^\infty}||n||_{L^\infty}+1)||m||_{H^s}
\end{eqnarray}
and
\begin{eqnarray}
||uv-u_x v_x||_{L^{\infty}}+||u_x n+v_x m||_{L^{\infty}}
\leq C ||m||_{L^\infty}||n||_{L^\infty}.
\end{eqnarray}
Hence,
\begin{eqnarray*}
||m(t)||_{H^{s}}\leq ||m_0||_{H^{s}}+C \int_0^t
(||m(\tau)||_{L^\infty}||n(\tau)||_{L^\infty}+1)||m(\tau)||_{H^s} d\tau.
\end{eqnarray*}
Similarly,
\begin{eqnarray*}
||n(t)||_{H^{s}}\leq ||n_0||_{H^{s}}+C \int_0^t
(||m(\tau)||_{L^\infty}||n(\tau)||_{L^\infty}+1)||n(\tau)||_{H^s} d\tau.
\end{eqnarray*}
Thus, we have
\begin{eqnarray}
\,\, && ||m(t)||_{H^{s}}+||n(t)||_{H^{s}}\\ \nonumber
&\leq& ||m_0||_{H^{s}}+||n_0||_{H^{s}}+C \int_0^t
(||m||_{L^\infty}||n||_{L^\infty}+1)(||m||_{H^s}+||n||_{H^s}) d\tau.
\end{eqnarray}
Taking advantage of Gronwall's inequality, one gets
\begin{eqnarray}
||m(t)||_{H^{s}}+||n(t)||_{H^{s}}
&\leq& (||m_0||_{H^{s}}+||n_0||_{H^{s}})\\ \nonumber
&&\times e^{C \int_0^t (||m||_{L^\infty}||n||_{L^\infty}+1) d\tau}.
\end{eqnarray}
Therefore, if $T<\infty$ satisfies $\int_0^T ||m(\tau)||_{L^\infty}||n(\tau)||_{L^\infty} d\tau<\infty$, then we deduce from (4.8) that
\begin{eqnarray}
\limsup\limits_{t\to T}(||m(t)||_{H^s}+||n(t)||_{H^s})<\infty,
\end{eqnarray}
which contradicts the assumption that $T<\infty$ is the maximal existence time. This completes the proof of the theorem for $s\in(\frac 1 2, 1)$.\\

\noindent {\it{Step 2}}. For $s\in[1,\frac{3}{2})$,
Lemma 2.1 (1) applied to the first equation of System (\ref{mn}), we get
\begin{eqnarray*}
||m(t)||_{H^{s}}&\leq& ||m_0||_{H^{s}}
+C\int_0^t ||m(\tau)||_{H^s}||u_x n+v_x m||_{H^{\frac 1 2}\cap L^{\infty}} d\tau\\
&&+C \int_0^t
(||\big((u_x n+v_x m)-(uv_x-u_xv)\big)m(\tau)||_{H^s}+||u_x(\tau)||_{H^s})d\tau.
\end{eqnarray*}
Noting that
\begin{eqnarray*}
||u_x n+v_x m||_{H^{\frac 1 2}\cap L^{\infty}}
\leq C ||u_x n+v_x m||_{H^{{\frac {1} {2}}+\varepsilon_0}}
\leq C ||m||_{H^{{\frac {1} {2}}+\varepsilon_0}} ||n||_{H^{{\frac {1} {2}}+\varepsilon_0}},
\end{eqnarray*}
where $\varepsilon_0\in (0,\frac 1 2)$ and we used the fact that
$H^{\frac 1 2 +\varepsilon_0}(\mathbb{R})\hookrightarrow H^{\frac 1 2}(\mathbb{R})\cap L^\infty(\mathbb{R})$,
which together with (4.5) yields
\begin{eqnarray*}
||m(t)||_{H^{s}}\leq ||m_0||_{H^{s}}+C \int_0^t
(||m||_{H^{{\frac {1} {2}}+\varepsilon_0}} ||n||_{H^{{\frac {1} {2}}+\varepsilon_0}}+1)||m(\tau)||_{H^s} d\tau.
\end{eqnarray*}
For the second equation of System (\ref{mn}), we can deal with it in a similar way and obtain that
\begin{eqnarray*}
||n(t)||_{H^{s}}\leq ||n_0||_{H^{s}}+C \int_0^t
(||m||_{H^{{\frac {1} {2}}+\varepsilon_0}} ||n||_{H^{{\frac {1} {2}}+\varepsilon_0}}+1)||n(\tau)||_{H^s} d\tau.
\end{eqnarray*}
Hence,
\begin{eqnarray*}
\quad  &&||m(t)||_{H^{s}}+||n(t)||_{H^{s}}\\ \nonumber
&\leq& ||m_0||_{H^{s}}+||n_0||_{H^{s}}+C \int_0^t
(||m||_{H^{{\frac {1} {2}}+\varepsilon_0}} ||n||_{H^{{\frac {1} {2}}+\varepsilon_0}}+1)
(||m||_{H^s}+||n||_{H^s}) d\tau.
\end{eqnarray*}
Thanks to Gronwall's inequality again, we have
\begin{eqnarray}
||m(t)||_{H^{s}}+||n(t)||_{H^{s}}
&\leq& (||m_0||_{H^{s}}+||n_0||_{H^{s}})\\ \nonumber
&&\times e^{C \int_0^t (||m||_{H^{{\frac {1} {2}}+\varepsilon_0}} ||n||_{H^{{\frac {1} {2}}+\varepsilon_0}}+1) d\tau}.
\end{eqnarray}
Therefore, if $T<\infty$ satisfies $\int_0^T ||m(\tau)||_{L^\infty}||n(\tau)||_{L^\infty} d\tau<\infty$, then we deduce from the uniqueness of the solution to System (\ref{mn}) and (4.9) with
$\frac{1}{2}+\varepsilon_0\in (\frac{1}{2},1)$ instead of $s$ that
$$||m(t)||_{H^{{\frac {1} {2}}+\varepsilon_0}} ||n(t)||_{H^{{\frac {1} {2}}+\varepsilon_0}}\quad
\text{is uniformly bounded in}\quad  t\in (0,T).$$
 This along with (4.10) implies that
\begin{eqnarray}
\limsup\limits_{t\to T}(||m(t)||_{H^s}+||n(t)||_{H^s})<\infty,
\end{eqnarray}
which contradicts the assumption that $T<\infty$ is the maximal existence time. This completes the proof of the theorem for $s\in[1,\frac 3 2)$.\\

\noindent {\it{Step 3}}. For $s\in(1,2)$,
by differentiating System (\ref{mn}) with respect to $x$, we have
\begin{eqnarray*}
&&\partial_t m_x+ \frac{1}{2}(uv-u_x v_x)\partial_x m_x \\ \nonumber
&=&[\frac{1}{2}(uv_x-u_x v)-(u_x n+v_x m)]m_x
-\frac{1}{2}[(u_x n+v_x m)-(uv_x-u_x v)]_x m-b u_{xx}\\ \nonumber
&\triangleq& R_1(t,x)
\end{eqnarray*}
and
\begin{eqnarray*}
&&\partial_t n_x+ \frac{1}{2}(uv-u_x v_x)\partial_x n_x \\ \nonumber
&=&[-\frac{1}{2}(uv_x-u_x v)-(u_x n+v_x m)]n_x
-\frac{1}{2}[(u_x n+v_x m)+(uv_x-u_x v)]_x n-b v_{xx}\\ \nonumber
&\triangleq& R_2(t,x).
\end{eqnarray*}
By Lemma 2.2 with $s-1\in(0,1)$, we get
\begin{eqnarray*}
||m_x(t)||_{H^{s-1}}&\leq& ||\partial_x m_0||_{H^{s-1}}
+C \int_0^t ||R_1(\tau)||_{H^{s-1}}d\tau\\
&+& C\int_0^t ||m_x(\tau)||_{H^{s-1}}(||uv-u_x v_x||_{L^{\infty}}+||u_x n+v_x m||_{L^{\infty}}) d\tau
\end{eqnarray*}
and
\begin{eqnarray*}
||n_x(t)||_{H^{s-1}}&\leq& ||\partial_x n_0||_{H^{s-1}}
+C \int_0^t ||R_2(\tau)||_{H^{s-1}}d\tau\\
&+& C\int_0^t ||n_x(\tau)||_{H^{s-1}}(||uv-u_x v_x||_{L^{\infty}}+||u_x n+v_x m||_{L^{\infty}}) d\tau.
\end{eqnarray*}
Thanks to (2.3) and (4.1)-(4.4), we have
\begin{eqnarray*}
&&||[\frac{1}{2}(uv_x-u_x v)-(u_x n+v_x m)]m_x||_{H^{s-1}}\\ \nonumber
&\leq& C(||\frac{1}{2}(uv_x-u_x v)-(u_x n+v_x m)||_{H^{s}}||m||_{L^{\infty}}\\ \nonumber
&&+||\frac{1}{2}(uv_x-u_x v)-(u_x n+v_x m)||_{L^{\infty}}||m_x||_{H^{s-1}})\\ \nonumber
&\leq& C ||m||_{L^{\infty}}||n||_{L^{\infty}}||m||_{H^{s}},
\end{eqnarray*}
and
\begin{eqnarray*}
&&||-\frac{1}{2}[(u_x n+v_x m)-(uv_x-u_x v)]_x m+ b u_{xx}||_{H^{s-1}}\\ \nonumber
&\leq& C(||m||_{H^{s}}||(u_x n+v_x m)-(uv_x-u_x v)||_{L^{\infty}}\\ \nonumber
&&+||m||_{L^{\infty}}||(u_x n+v_x m)-(uv_x-u_x v)||_{H^{s}})+|b|||m||_{H^{s-1}}\\ \nonumber
&\leq& C (||m||_{L^{\infty}}||n||_{L^{\infty}}+1)||m||_{H^{s}},
\end{eqnarray*}
which together with (4.6) yields
\begin{eqnarray*}
||m_x(t)||_{H^{s-1}}\leq ||m_0||_{H^{s}}+C \int_0^t
(||m(\tau)||_{L^\infty}||n(\tau)||_{L^\infty}+1)||m(\tau)||_{H^s} d\tau.
\end{eqnarray*}
Similarly,
\begin{eqnarray*}
||n_x(t)||_{H^{s-1}}\leq ||n_0||_{H^{s}}+C \int_0^t
(||m(\tau)||_{L^\infty}||n(\tau)||_{L^\infty}+1)||n(\tau)||_{H^s} d\tau.
\end{eqnarray*}
Thus, we have
\begin{eqnarray*}
\,\, && ||m_x(t)||_{H^{s-1}}+||n_x(t)||_{H^{s-1}}\\ \nonumber
&\leq& ||m_0||_{H^{s}}+||n_0||_{H^{s}}+C \int_0^t
(||m||_{L^\infty}||n||_{L^\infty}+1)(||m||_{H^s}+||n||_{H^s}) d\tau.
\end{eqnarray*}
This along with (4.7) with $s-1$ instead of $s$ ensures
\begin{eqnarray*}
\,\, && ||m(t)||_{H^{s}}+||n(t)||_{H^{s}}\\ \nonumber
&\leq& ||m_0||_{H^{s}}+||n_0||_{H^{s}}+C \int_0^t
(||m||_{L^\infty}||n||_{L^\infty}+1)(||m||_{H^s}+||n||_{H^s}) d\tau.
\end{eqnarray*}
Similar to Step 1, we can easily prove the theorem for $s\in(1,2)$.\\

\noindent {\it{Step 4}}. For $s=k\in\mathbb{N}$ and $k\geq2$,
by differentiating System (\ref{mn}) $k-1$ times with respect to $x$, we get
\begin{eqnarray*}
(\partial_t+\frac{1}{2}(uv-u_x v_x)\partial_x)\partial^{k-1}_x m
&=&-\frac{1}{2}\sum\limits_{l=0}^{k-2} C^l_{k-1}\partial^{k-l-1}_x(uv-u_x v_x)\partial^{l+1}_x m
-b \partial^{k-1}_x u_x\\ \nonumber
&&-\frac{1}{2}\partial^{k-1}_x ([(u_x n+v_x m)-(uv_x-u_x v)]m)\\ \nonumber
&\triangleq& F_1(t,x)
\end{eqnarray*}
and
\begin{eqnarray*}
(\partial_t+\frac{1}{2}(uv-u_x v_x)\partial_x)\partial^{k-1}_x n
&=&-\frac{1}{2}\sum\limits_{l=0}^{k-2} C^l_{k-1}\partial^{k-l-1}_x(uv-u_x v_x)\partial^{l+1}_x n
-b \partial^{k-1}_x v_x\\ \nonumber
&&-\frac{1}{2}\partial^{k-1}_x ([(u_x n+v_x m)+(uv_x-u_x v)]n)\\ \nonumber
&\triangleq& F_2(t,x),
\end{eqnarray*}
which together with Lemma 2.1 (1), implies that
\begin{eqnarray*}
||\partial^{k-1}_x m(t)||_{H^1}
\leq ||m_0||_{H^k}+ \int_0^t ||F_1(\tau)||_{H^1}d\tau
+C \int_0^t ||u_x n+v_x m||_{H^{\frac 1 2}\cap L^{\infty}}||m||_{H^k} d\tau
\end{eqnarray*}
and
\begin{eqnarray*}
||\partial^{k-1}_x n(t)||_{H^1}
\leq ||n_0||_{H^k}+ \int_0^t ||F_2(\tau)||_{H^1}d\tau
+C \int_0^t ||u_x n+v_x m||_{H^{\frac 1 2}\cap L^{\infty}}||n||_{H^k} d\tau.
\end{eqnarray*}
Thanks to (4.1)-(4.4) again, we have
\begin{eqnarray}
&&||-\frac{1}{2}\sum\limits_{l=0}^{k-2} C^l_{k-1}\partial^{k-l-1}_x(uv-u_x v_x)\partial^{l+1}_x m||_{H^{1}}
\\ \nonumber
&\leq& C(k)\sum\limits_{l=0}^{k-2}||\partial^{k-l-1}_x(uv-u_x v_x)||_{L^{\infty}}||m||_{H^{l+2}}\\ \nonumber
&\leq& C(k)\sum\limits_{l=0}^{k-2}||uv-u_x v_x||_{H^{k-l-\frac{1}{2}+\varepsilon_0}}||m||_{H^{l+2}}\\ \nonumber
&\leq& C(k) ||uv-u_x v_x||_{H^{k-\frac{1}{2}+\varepsilon_0}}||m||_{H^{k}} \\ \nonumber
&\leq& C(k) ||m||_{H^{k-\frac{1}{2}+\varepsilon_0}}||n||_{H^{k-\frac{1}{2}+\varepsilon_0}}||m||_{H^{k}},
\end{eqnarray}
\begin{eqnarray}
&&||-\frac{1}{2}\partial^{k-1}_x ([(u_x n+v_x m)-(uv_x-u_x v)]m)-b \partial^{k-1}_x u_x||_{H^{1}}\\ \nonumber
&\leq& C||[(u_x n+v_x m)-(uv_x-u_x v)]m||_{H^{k}}+ |b|||u_x||_{H^{k}}\\ \nonumber
&\leq& C(||m||_{L^{\infty}}||m||_{L^{\infty}}+1)||m||_{H^{k}}\\ \nonumber
&\leq& C(||m||_{H^{k-\frac{1}{2}+\varepsilon_0}}||n||_{H^{k-\frac{1}{2}+\varepsilon_0}}+1)||m||_{H^{k}},
\end{eqnarray}
and
\begin{eqnarray*}
||u_x n+v_x m||_{H^{\frac{1}{2}}\cap L^{\infty}}
\leq C ||u_x n+v_x m||_{H^{k-\frac{1}{2}+\varepsilon_0}}
\leq C ||m||_{H^{k-\frac{1}{2}+\varepsilon_0}}||n||_{H^{k-\frac{1}{2}+\varepsilon_0}},
\end{eqnarray*}
where $\varepsilon_0\in (0,\frac 1 2)$ and we used the fact that
\begin{eqnarray}
H^{k-\frac 1 2 +\varepsilon_0}(\mathbb{R})\hookrightarrow H^{\frac 1 2 +\varepsilon_0}(\mathbb{R})\hookrightarrow H^{\frac 1 2}(\mathbb{R})\cap L^\infty(\mathbb{R})\quad \text{with}\quad  k\geq 2.
\end{eqnarray}
Thus, we get
\begin{eqnarray*}
||\partial^{k-1}_x m(t)||_{H^{1}}\leq ||m_0||_{H^{k}}+C \int_0^t
(||m||_{H^{k-\frac{1}{2}+\varepsilon_0}}||n||_{H^{k-\frac{1}{2}+\varepsilon_0}}+1)||m||_{H^k} d\tau.
\end{eqnarray*}
Similarly,
\begin{eqnarray*}
||\partial^{k-1}_x n(t)||_{H^{1}}\leq ||n_0||_{H^{k}}+C \int_0^t
(||m||_{H^{k-\frac{1}{2}+\varepsilon_0}}||n||_{H^{k-\frac{1}{2}+\varepsilon_0}}+1)||n||_{H^k} d\tau.
\end{eqnarray*}
Then,
\begin{eqnarray*}
\,\, && ||\partial^{k-1}_x m(t)||_{H^{1}}+||\partial^{k-1}_x n(t)||_{H^{1}}\\ \nonumber
&\leq& ||m_0||_{H^{k}}+||n_0||_{H^{k}}+C \int_0^t
(||m||_{H^{k-\frac{1}{2}+\varepsilon_0}}||n||_{H^{k-\frac{1}{2}+\varepsilon_0}}+1)(||m||_{H^k}+||n||_{H^k}) d\tau,
\end{eqnarray*}
which together with Gronwall's inequality and (4.10) with $s=1$, yields that
\begin{eqnarray}
||m(t)||_{H^{k}}+||n(t)||_{H^{k}}
&\leq& (||m_0||_{H^{k}}+||n_0||_{H^{k}})\\ \nonumber
&&\times e^{C \int_0^t
(||m||_{H^{k-\frac{1}{2}+\varepsilon_0}}||n||_{H^{k-\frac{1}{2}+\varepsilon_0}}+1) d\tau}.
\end{eqnarray}
If $T<\infty$ satisfies $\int_0^T ||m(\tau)||_{L^\infty}||n(\tau)||_{L^\infty} d\tau<\infty$,
applying Step 3 with $\frac 3 2 +\varepsilon_0\in(1,2)$ and by induction with respect to $k\geq 2$, we can obtain that $||m(t)||_{H^{k-\frac{1}{2}+\varepsilon_0}}||n(t)||_{H^{k-\frac{1}{2}+\varepsilon_0}}$
is uniformly bounded in $t\in (0,T)$.
Thanks to (4.15), we get
\begin{eqnarray}
\limsup\limits_{t\to T}(||m(t)||_{H^k}+||n(t)||_{H^k})<\infty,
\end{eqnarray}
which contradicts the assumption that $T<\infty$ is the maximal existence time. This completes the proof of the theorem for $s=k\in\mathbb{N}$ and $k\geq 2$.\\

\noindent {\it{Step 5}}. For $s\in(k,k+1),\,k\in\mathbb{N} $ and $k\geq 2$,
by differentiating System (\ref{mn}) $k$ times with respect to $x$, we get
\begin{eqnarray*}
(\partial_t+\frac{1}{2}(uv-u_x v_x)\partial_x)\partial^{k}_x m
&=&-\frac{1}{2}\sum\limits_{l=0}^{k-1} C^l_{k}\partial^{k-l}_x(uv-u_x v_x)\partial^{l+1}_x m
-b \partial^{k}_x u_x\\ \nonumber
&&-\frac{1}{2}\partial^{k}_x ([(u_x n+v_x m)-(uv_x-u_x v)]m)\\ \nonumber
&\triangleq& G_1(t,x)
\end{eqnarray*}
and
\begin{eqnarray*}
(\partial_t+\frac{1}{2}(uv-u_x v_x)\partial_x)\partial^{k}_x n
&=&-\frac{1}{2}\sum\limits_{l=0}^{k-1} C^l_{k}\partial^{k-l}_x(uv-u_x v_x)\partial^{l+1}_x n
-b \partial^{k}_x v_x\\ \nonumber
&&-\frac{1}{2}\partial^{k}_x ([(u_x n+v_x m)+(uv_x-u_x v)]n)\\ \nonumber
&\triangleq& G_2(t,x),
\end{eqnarray*}
which together with Lemma 2.2 with $s-k\in(0,1)$, implies that
\begin{eqnarray*}
||\partial^{k}_x m(t)||_{H^{s-k}}
&\leq& ||\partial^{k}_x m_0||_{H^{s-k}}+ C \int_0^t ||G_1(\tau)||_{H^{s-k}}d\tau\\ \nonumber
&&+C \int_0^t (||u v-u_x v_x||_{L^{\infty}}+||u_x n+v_x m||_{L^{\infty}})||\partial^{k}_x m(\tau)||_{H^{s-k}} d\tau
\end{eqnarray*}
and
\begin{eqnarray*}
||\partial^{k}_x n(t)||_{H^{s-k}}
&\leq& ||\partial^{k}_x n_0||_{H^{s-k}}+ C \int_0^t ||G_2(\tau)||_{H^{s-k}}d\tau\\ \nonumber
&&+C \int_0^t (||u v-u_x v_x||_{L^{\infty}}+||u_x n+v_x m||_{L^{\infty}})||\partial^{k}_x n(\tau)||_{H^{s-k}} d\tau.
\end{eqnarray*}
Thanks to (4.14) again and similar to (4.12)-(4.13), we have
\begin{eqnarray*}
&&||-\frac{1}{2}\partial^{k}_x ([(u_x n+v_x m)-(uv_x-u_x v)]m)-b \partial^{k}_x u_x||_{H^{s-k}}\\ \nonumber
&\leq& C(||m||_{H^{k-\frac{1}{2}+\varepsilon_0}}||n||_{H^{k-\frac{1}{2}+\varepsilon_0}}+1)||m||_{H^{s}},
\end{eqnarray*}
and
\begin{eqnarray*}
||-\frac{1}{2}\sum\limits_{l=1}^{k-1} C^l_{k}\partial^{k-l}_x(uv-u_x v_x)\partial^{l+1}_x m||_{H^{s-k}}
\leq C(k) ||m||_{H^{k-\frac{1}{2}+\varepsilon_0}}||n||_{H^{k-\frac{1}{2}+\varepsilon_0}}||m||_{H^{s}}.
\end{eqnarray*}
Applying (2.3), (4.1)-(4.4) and (4.14), one infers
\begin{eqnarray*}
&&||-\frac{1}{2} C^0_{k}\partial^{k}_x(uv-u_x v_x) m_x||_{H^{s-k}}\\ \nonumber
&\leq& C (||m_x||_{H^{s-k+1}}||\partial^{k-1}_x(uv-u_x v_x)||_{L^\infty}
+||m_x||_{L^\infty}||\partial^{k}_x(uv-u_x v_x)||_{H^{s-k}})\\ \nonumber
&\leq& C (||m||_{H^{s-k+2}}||uv-u_x v_x||_{H^{k-\frac{1}{2}+\varepsilon_0}}
+||m||_{H^{k-\frac{1}{2}+\varepsilon_0}}||uv-u_x v_x||_{H^{s}})\\ \nonumber
&\leq& C ||m||_{H^{k-\frac{1}{2}+\varepsilon_0}}||n||_{H^{k-\frac{1}{2}+\varepsilon_0}}||m||_{H^{s}}.
\end{eqnarray*}
Thus, we get
\begin{eqnarray*}
||\partial^{k}_x m(t)||_{H^{s-k}}\leq ||m_0||_{H^{s}}+C \int_0^t
(||m||_{H^{k-\frac{1}{2}+\varepsilon_0}}||n||_{H^{k-\frac{1}{2}+\varepsilon_0}}+1)||m||_{H^s} d\tau.
\end{eqnarray*}
Similarly,
\begin{eqnarray*}
||\partial^{k}_x n(t)||_{H^{s-k}}\leq ||n_0||_{H^{s}}+C \int_0^t
(||m||_{H^{k-\frac{1}{2}+\varepsilon_0}}||n||_{H^{k-\frac{1}{2}+\varepsilon_0}}+1)||n||_{H^s} d\tau.
\end{eqnarray*}
Then,
\begin{eqnarray*}
\,\, && ||\partial^{k}_x m(t)||_{H^{s-k}}+||\partial^{k}_x n(t)||_{H^{s-k}}\\ \nonumber
&\leq& ||m_0||_{H^{s}}+||n_0||_{H^{s}}+C \int_0^t
(||m||_{H^{k-\frac{1}{2}+\varepsilon_0}}||n||_{H^{k-\frac{1}{2}+\varepsilon_0}}+1)(||m||_{H^s}+||n||_{H^s}) d\tau,
\end{eqnarray*}
which along with (4.7) with $s-k\in (0,1)$ instead of $s$, ensures that
\begin{eqnarray*}
\,\, && ||m(t)||_{H^{s}}+|| n(t)||_{H^{s}}\\ \nonumber
&\leq& ||m_0||_{H^{s}}+||n_0||_{H^{s}}+C \int_0^t
(||m||_{H^{k-\frac{1}{2}+\varepsilon_0}}||n||_{H^{k-\frac{1}{2}+\varepsilon_0}}+1)(||m||_{H^s}+||n||_{H^s}) d\tau,
\end{eqnarray*}
By using Gronwall's inequality, Step 3 with $\frac 3 2+\varepsilon_0\in (1,2)$ and the similar argument as that in Step 4, one can easily get the desired result.

Consequently, we have completed the proof of the theorem from Step 1 to Step 5.
\end{proof}

\begin{remark4}
The maximal existence time $T$ in Theorem 4.1 can be chosen independent of the regularity index $s$. Indeed, let
$(m_0,n_0)\in H^{s}\times H^{s}$ with $s>\frac{1}{2}$ and some $s'\in (\frac 1 2,s)$. Then Remark 3.1 ensures that there exists a unique $H^s\times H^s$ (resp., $H^{s'}\times H^{s'}$) solution $(m_s,n_s)$ (resp., $(m_{s'},n_{s'})$) to System (\ref{mn}) with the maximal existence time $T_s$ (resp., $T_{s'}$). Since $H^s\hookrightarrow H^{s'}$, it follows from the uniqueness that $T_s\leq T_{s'}$ and $(m_s,n_s)\equiv (m_{s'},n_{s'})$ on $[0, T_s)$. On the other hand, if we suppose that $T_s < T_{s'}$, then $(m_{s'},n_{s'})\in C([0,T_s]; H^{s'}\times H^{s'})$.
Hence $(m_{s'},n_{s'})\in L^2(0,T_s; L^{\infty}\times L^{\infty})$, which together with the Holder inequality leads to a contraction to Theorem 4.1. Therefore, $T_s = T_{s'}$.
\end{remark4}

Making use of Sobolev's embedding theorem and Theorem 4.1, one can readily get the blow-up criterion as follows.
\begin{corollary4}
Let $(m_0,n_0)\in H^s(\mathbb{R})\times H^{s}(\mathbb{R})$
with $s> \frac 1 2$ and $T>0$ be the maximal existence time of the corresponding solution
$(m,n)$ to System (\ref{mn}), which is guaranteed by Remark 3.1.
Then the solution $(m,n)$ blows up in finite time
if and only if
$$\limsup_{t\rightarrow T}||m(t,\cdot)||_{L^\infty}=\infty \quad \text{or}\quad
\limsup_{t\rightarrow T}||n(t,\cdot)||_{L^\infty}=\infty.$$
\end{corollary4}

\par
Now we turn our attention to the precise blow-up scenario for sufficiently regular solutions to System (\ref{mn}) with $b=0$ which will be assumed for the remainder of this section. For this, we firstly consider the following ordinary differential equation:
\begin{equation}
\left\{\begin{array}{ll}
\frac{d q(t,x)}{dt}=\frac{1}{2}(u v-u_x v_x)(t,q(t,x)),\quad\quad &(t,x)\in(0,T)\times \mathbb{R}, \\
q(0,x)=x,\quad\quad  &x\in\mathbb{R},
\end{array}\right.
\end{equation}
for the flow generated by $\frac{1}{2}(u v-u_x v_x)$.

The following lemmas are very crucial to study the blow-up phenomena of strong solutions to System (\ref{mn}) with $b=0$.

\begin{lemma4}
Let $(m_0,n_0)\in H^s(\mathbb{R})\times H^s(\mathbb{R})$ with $s> \frac 1 2$
and $T>0$ be the maximal existence time of the corresponding solution $(m,n)$ to System (\ref{mn}) with $b=0$.
Then Eq.(4.17) has a unique solution $q\in C^1([0,T)\times \mathbb{R};\mathbb{R})$.
Moreover, the mapping $q(t,\cdot)$ is an increasing diffeomorphism of $\mathbb{R}$ with
\begin{eqnarray}
q_{x}(t,x)=\exp\left(\frac{1}{2}\int_0^t (u_x n+v_x m)(s,q(s,x))ds\right)>0,
\end{eqnarray}
for all $(t,x)\in [0,T)\times \mathbb{R}$.
\end{lemma4}

\begin{proof}
Since $(u,v)\in C([0,T); H^s(\mathbb{R})\times H^s(\mathbb{R}))
\cap C^1([0,T); H^{s-1}(\mathbb{R})\times H^{s-1}(\mathbb{R}))$ with $s>\frac 5 2$, it follows from the fact
$H^{s-1}(\mathbb{R})\hookrightarrow Lip(\mathbb{R})$ with $s>\frac 5 2$ that $\frac{1}{2}(u v-u_x v_x)$ is bounded and Lipschitz continuous in the space variable $x$ and of class $C^1$ in time variable $t$. Then the classical ODE theory ensures that Eq.(4.17) has a unique solution $q\in C^1([0,T)\times \mathbb{R};\mathbb{R})$.\\
By differentiating Eq.(4.17) with respect to $x$, one gets
\begin{equation*}
\left\{\begin{array}{ll}
\frac{d q_x(t,x)}{dt}=\frac{1}{2}(u_x n+v_x m)(t,q(t,x))q_x(t,x),\quad\quad &(t,x)\in(0,T)\times \mathbb{R}, \\
q_x(0,x)=1,\quad\quad  &x\in\mathbb{R},
\end{array}\right.
\end{equation*}
which leads to (4.18).\\
On the other hand, $\forall\, t<T$, by the Sobolev embedding theorem, we have
$$\sup\limits_{(s,x)\in [0,T)\times \mathbb{R}} |\frac{1}{2}(u_x n+v_x m)(s,x)| < \infty,$$
which along with (4.18) yields that there exists a constant $C>0$ such that
$$q_x(t,x)\geq e^{-Ct},\quad \forall\, (t,x)\in [0,T)\times \mathbb{R}.$$
This implies that  the mapping $q(t,\cdot)$ is an increasing diffeomorphism of $\mathbb{R}$ before blow-up.
Therefore, we complete the proof of Lemma 4.1.
\end{proof}

\begin{lemma4}
Let $(m_0,n_0)\in H^s(\mathbb{R})\times H^s(\mathbb{R})$ with $s> \frac 1 2$
and $T>0$ be the maximal existence time of the corresponding solution $(m,n)$ to System (\ref{mn}) with $b=0$.
Then we have
\begin{equation}
m(t,q(t,x))q_x(t,x)=m_0(x) \exp\left(\frac{1}{2}\int_0^t (u v_x-u_x v)(s,q(s,x))ds\right),
\end{equation}
and
\begin{equation}
n(t,q(t,x))q_x(t,x)=n_0(x) \exp\left(-\frac{1}{2}\int_0^t (u v_x-u_x v)(s,q(s,x))ds\right).
\end{equation}
for all $(t,x)\in [0,T)\times \mathbb{R}$.\\
Moreover, if there exists an $C>0$ such that $(u_x n+v_x m)(t,x)\geq -C$
and $||(u v_x-u_x v)(t,\cdot)||_{L^\infty}\leq C$ for all $(t,x)\in [0,T)\times \mathbb{R}$, then
\begin{equation}
||m(t,\cdot)||_{L^\infty}\leq C e^{Ct}||m_0||_{H^s} \quad \text{and} \quad
||n(t,\cdot)||_{L^\infty}\leq C e^{Ct}||n_0||_{H^s},
\end{equation}
for all $t\in [0,T)$.
\end{lemma4}

\begin{proof}
By differentiating the left-hand side of (4.19)-(4.20) with
respect to $t$ and making use of (4.17)-(4.18) and System (\ref{mn}), one infers
\begin{eqnarray*}
&&\frac{d}{dt}{(m(t,q(t,x))q_x(t,x))}\\
&=&(m_t(t,q)+m_x(t,q)q_t(t,x))q_{x}(t,x)+m(t,q)q_{xt}(t,x)\\
&=&\big(m_t+\frac{1}{2}(u v-u_x v_x) m_x+\frac{1}{2}(u_x n+v_x m) m\big)(t,q(t,x)) q_{x}(t,x)\\
&=& \frac{1}{2}(u v_x-u_x v)(t,q(t,x))m(t,q(t,x))q_x(t,x)
\end{eqnarray*}
and
\begin{eqnarray*}
&&\frac{d}{dt}{(n(t,q(t,x))q_x(t,x))}\\
&=&(n_t(t,q)+n_x(t,q)q_t(t,x))q_{x}(t,x)+n(t,q)q_{xt}(t,x)\\
&=&\big(n_t+\frac{1}{2}(u v-u_x v_x) n_x+\frac{1}{2}(u_x n+v_x m) n\big)(t,q(t,x)) q_{x}(t,x)\\
&=& -\frac{1}{2}(u v_x-u_x v)(t,q(t,x))n(t,q(t,x))q_x(t,x),
\end{eqnarray*}
which proves (4.19) and (4.20). By Lemma 4.1, in view of (4.18)-(4.20) and the assumption of the lemma,
we obtain for all $t\in[0,T)$,
\begin{eqnarray*}
||m(t,\cdot)||_{L^{\infty}}
&=&||m(t,q(t,\cdot))||_{L^{\infty}}\\
&=&||e^{\frac{1}{2}\int_0^t (u v_x-u_x v)(s,\cdot)ds}q_x^{-1}(t,\cdot) m_0(\cdot)||_{L^{\infty}}\\
&\leq& C e^{Ct}||m_0||_{H^{s}}
\end{eqnarray*}
and
\begin{eqnarray*}
||n(t,\cdot)||_{L^{\infty}}
&=&||n(t,q(t,\cdot))||_{L^{\infty}}\\
&=&||e^{-\frac{1}{2}\int_0^t (u v_x-u_x v)(s,\cdot)ds}q_x^{-1}(t,\cdot) n_0(\cdot)||_{L^{\infty}}\\
&\leq& C e^{Ct}||n_0||_{H^{s}}.
\end{eqnarray*}
Thus, we completes the proof of the lemma.
\end{proof}

The following theorem shows the precise blow-up scenario for sufficiently regular solutions
to System (\ref{mn}) with $b=0$.

\begin{theorem4}
Let $(m_0,n_0)\in H^s(\mathbb{R})\times H^s(\mathbb{R})$ with $s> \frac 1 2$
and $T>0$ be the maximal existence time of the corresponding solution $(m,n)$ to System (\ref{mn}) with $b=0$.
Then the solution $(m,n)$ blows up in finite time
if and only if
$$\liminf_{t\rightarrow T}\inf_{x\in\mathbb{R}}\{(u_x n+ v_x m)(t,x)\}=-\infty \quad \text{or}\quad
\limsup_{t\rightarrow T} (||(u v_x-u_x v)(t,\cdot)||_{L^\infty})=\infty.$$
\end{theorem4}

\begin{proof}
Assume that the solution $(m,n)$ blows up in finite time ($T<\infty$) and
there exists a $C>0$ such that
$$(u_x n+v_x m)(t,x)\geq -C\quad \text{and} \quad
||(u v_x-u_x v)(t,\cdot)||_{L^\infty}\leq C,\quad \forall\, (t,x)\in [0,T)\times \mathbb{R}.$$
By (4.21), we have
$$\int_0^T ||m(t)||_{L^\infty}||n(t)||_{L^\infty} d t\leq C^2 T e^{2CT}||m_0||_{H^{s}}||n_0||_{H^{s}}<\infty,$$
which contradicts to Theorem 4.1.

On the other hand, by Sobolev's embedding theorem, we can see that if
$$\liminf_{t\rightarrow T}\inf_{x\in\mathbb{R}}\{(u_x n+ v_x m)(t,x)\}=-\infty \quad \text{or}\quad
\limsup_{t\rightarrow T} (||(u v_x-u_x v)(t,\cdot)||_{L^\infty})=\infty,$$
then the solution $(m,n)$ will blow up in finite time. This completes the proof of the theorem.
\end{proof}

\begin{remark4}
If $v=2u$, then Theorem 4.2 covers the corresponding result in \cite{Gui-CMP}.
\end{remark4}

In order to state a new blow-up criterion with respect to the initial data of strong solutions to System (\ref{mn}) with $b=0$, we first investigate the transport equation in terms of $\frac{1}{2}(u_x n+v_x m)$ which is the slope of
$\frac{1}{2}(u v-u_x v_x)$.

\begin{lemma4}
Let $(m_0,n_0)\in H^s(\mathbb{R})\times H^s(\mathbb{R})$ with $s> \frac 1 2$
and $T>0$ be the maximal existence time of the corresponding solution $(m,n)$ to System (\ref{mn}) with $b=0$.
Set $M=M(t,x)\triangleq (u_x n+v_x m)(t,x)$. Then for all $(t,x)\in [0,T)\times \mathbb{R}$, we have
\begin{eqnarray}
&&\quad\quad\quad  M_t+\frac{1}{2}(u v-u_x v_x) M_x\\ \nonumber
&=&-\frac{1}{2} M^2
-\frac{1}{2}n(1-\partial^2_x)^{-1}(u_x M)-\frac{1}{2}m(1-\partial^2_x)^{-1}(v_x M)
-\frac{1}{2}n\partial_x(1-\partial^2_x)^{-1}(u M)\\ \nonumber
&&-\frac{1}{2}m\partial_x(1-\partial^2_x)^{-1}(v M)
-\frac{1}{2}(u v_x-u_x v)(u_x n-v_x m)\\ \nonumber
&&+\frac{1}{2}n\partial_x(1-\partial^2_x)^{-1}((u v_x-u_x v)m)
-\frac{1}{2}m\partial_x(1-\partial^2_x)^{-1}((u v_x-u_x v)n).
\end{eqnarray}
Moreover, if we assume that $m_0(x),\, n_0(x)\geq 0$ for all $x\in\mathbb{R}$, then
$$|u_x(t,x)|\leq u(t,x),\quad   |v_x(t,x)|\leq v(t,x)   \quad \text{and}$$
\begin{eqnarray*}
 M_t+\frac{1}{2}(u v-u_x v_x) M_x
\leq -\frac{1}{2} M^2+\frac{7}{2}||u||_{L^\infty}||v||^2_{L^\infty}m+\frac{7}{2}||u||^2_{L^\infty}||v||_{L^\infty}n,
\end{eqnarray*}
for all $(t,x)\in [0,T)\times \mathbb{R}$.
\end{lemma4}

\begin{proof}
In view of Remark 4.1, we here may assume $s\geq 3$ to prove the lemma. Firstly, we have
\begin{eqnarray}
&&\quad\quad M_t+\frac{1}{2}(u v-u_x v_x) M_x \\ \nonumber
&=&u_{xt}n+v_{xt}m+u_x n_t+v_x m_t+\frac{1}{2}(u v-u_x v_x)(u_x n_x+v_x m_x+u_{xx}n+v_{xx}m).
\end{eqnarray}
From System (\ref{mn}), we infer that
\begin{eqnarray*}
&&(1-\partial^2_x)(u_t+\frac{1}{2}(u v-u_x v_x) u_x)\\
&=& m_t+\frac{1}{2}(1-\partial^2_x)((u v-u_x v_x) u_x)\\
&=& -\frac{1}{2}(u v-u_x v_x) m_x-\frac{1}{2}(u_x n+v_x m)m+\frac{1}{2}(u v_x-u_x v)m
+\frac{1}{2}(u v-u_x v_x) u_x\\
&&-\frac{1}{2}\partial^2_x ((u v-u_x v_x) u_x)\\
&=& -\frac{1}{2}(u_x n+v_x m)m+\frac{1}{2}(u v_x-u_x v)m-\frac{1}{2}(u_x n+v_x m)_x u_x-(u_x n+v_x m)u_{xx}\\
&=& -\frac{1}{2}(u M-(u v_x-u_x v)m+(u_x M)_x).
\end{eqnarray*}
Hence,
\begin{eqnarray}
&& (u_t+\frac{1}{2}(u v-u_x v_x) u_x)\\ \nonumber
&=&-\frac{1}{2}(1-\partial^2_x)^{-1}(u M-(u v_x-u_x v)m+(u_x M)_x).
\end{eqnarray}
Similarly,
\begin{eqnarray}
&& (v_t+\frac{1}{2}(u v-u_x v_x) v_x)\\ \nonumber
&=&-\frac{1}{2}(1-\partial^2_x)^{-1}(v M+(u v_x-u_x v)n+(v_x M)_x).
\end{eqnarray}
By virtue of (4.24)-(4.25) and System (\ref{mn}), one gets
\begin{eqnarray}
&& u_{xt}n+v_{xt}m \\ \nonumber
&=&-\frac{1}{2}(u v-u_x v_x)(u_{xx}n+v_{xx}m)
-\frac{1}{2}n\partial_x(1-\partial^2_x)^{-1}(u M-(u v_x-u_x v)m)\\ \nonumber
&&-\frac{1}{2}m\partial_x(1-\partial^2_x)^{-1}(v M+(u v_x-u_x v)n)
-\frac{1}{2}n(1-\partial^2_x)^{-1}(u_x M)\\ \nonumber
&&-\frac{1}{2}m(1-\partial^2_x)^{-1}(v_x M)
\end{eqnarray}
and
\begin{eqnarray*}
u_x n_t+v_x m_t
&=&-\frac{1}{2}(u v-u_x v_x)(u_x n_x+v_x m_x) \\ \nonumber
&&-\frac{1}{2}M^2-\frac{1}{2}(u v_x-u_x v)(u_x n-v_x m),
\end{eqnarray*}
which together with (4.23) and (4.26) reaches (4.22).\\
Since $m_0(x),\, n_0(x)\geq 0$ for all $x\in\mathbb{R}$, it follows from (4.18)-(4.20) that
\begin{eqnarray}
m(t,x),\, n(t,x)\geq 0, \quad \forall \, (t,x)\in[0,T)\times \mathbb{R}.
\end{eqnarray}
Note that
$$u(t,x)=(1-\partial^2_x)^{-1}m(t,x)=(p\ast m)(t,x)=\frac 1 2\int_{\mathbb{R}} e^{-|x-y|}m(t,y)d y.$$
Then
$$u(t,x)=\frac{e^{-x}}{2}\int^x_{-\infty} e^y m(t,y)d y+\frac{e^{x}}{2}\int_x^{\infty} e^{-y} m(t,y)d y
\quad \text{and}$$
$$u_x(t,x)=-\frac{e^{x}}{2}\int^x_{-\infty} e^y m(t,y)d y+\frac{e^{x}}{2}\int_x^{\infty} e^{-y} m(t,y)d y,$$
which together with (4.27) yields
$$u(t,x)+u_x(t,x)=e^x \int_x^{\infty} e^{-y} m(t,y)d y\geq 0\quad \quad \text{and}$$
$$u(t,x)-u_x(t,x)=e^{-x} \int^x_{-\infty} e^y m(t,y)d y\geq 0.$$
Hence,
\begin{eqnarray}
|u_x(t,x)|\leq u(t,x),\quad \forall \, (t,x)\in[0,T)\times \mathbb{R}.
\end{eqnarray}
Similarly,
\begin{eqnarray}
|v_x(t,x)|\leq v(t,x),\quad \forall \, (t,x)\in[0,T)\times \mathbb{R}.
\end{eqnarray}

Noting that
\begin{eqnarray*}
|\partial_x (1-\partial^2_x)^{-1}f(x)|
&=& |\frac 1 2 \int_{\mathbb{R}}sgn (x-y) e^{-|x-y|} f(y)d y|\\
&\leq& \frac 1 2 \int_{\mathbb{R}} e^{-|x-y|} |f(y)|d y\\
&=& (p\ast |f|)(x),
\end{eqnarray*}
and applying (4.28)-(4.29) and the facts that $u=p\ast m,\, v=p\ast n$, one obtains
\begin{eqnarray*}
&& -\frac{1}{2}n(1-\partial^2_x)^{-1}(u_x M)-\frac{1}{2}m(1-\partial^2_x)^{-1}(v_x M)\\
&\leq& -\frac{1}{2}n (p\ast (u_x v_x m))-\frac{1}{2}m (p\ast (u_x v_x n))\\
&\leq& \frac{1}{2}||u v||_{L^\infty}(u n+v m),
\end{eqnarray*}

\begin{eqnarray*}
&& -\frac{1}{2}n\partial_x(1-\partial^2_x)^{-1}(u M)-\frac{1}{2}m\partial_x(1-\partial^2_x)^{-1}(v M)\\
&\leq& \frac{1}{2}n
(||u^2||_{L^\infty}(p\ast n)+||u v||_{L^\infty}(p\ast m))
+\frac{1}{2}m
(||u v||_{L^\infty}(p\ast n)+||v^2||_{L^\infty}(p\ast m))\\
&=& \frac{1}{2}||u^2||_{L^\infty} v n+\frac{1}{2}||v^2||_{L^\infty} u m+\frac{1}{2}||u v||_{L^\infty}(u n+v m),
\end{eqnarray*}

\begin{eqnarray*}
-\frac{1}{2}(u v_x-u_x v)(u_x n-v_x m) \leq ||u v||_{L^\infty}(u n+v m),
\end{eqnarray*}
and
\begin{eqnarray*}
&& \frac{1}{2}n\partial_x(1-\partial^2_x)^{-1}((u v_x-u_x v)m)
-\frac{1}{2}m\partial_x(1-\partial^2_x)^{-1}((u v_x-u_x v)n)\\
&\leq& ||u v||_{L^\infty}(n (p\ast m)+m (p\ast n))\\
&=& ||u v||_{L^\infty}(u n+v m),
\end{eqnarray*}
which along with (4.22) implies the desired result. Therefore, we complete the proof of the lemma.
\end{proof}

Note that Theorem 4.2, Lemma 4.2 and (4.1)-(4.3) imply that if we want to study the fine structure of finite time singularities, one should assume in the following that there exists a constant $C=C(||m_0||_{H^s},||n_0||_{H^s})>0$ such that $||u(t,\cdot)||_{L^\infty},\,||v(t,\cdot)||_{L^\infty}\leq Ce^{Ct}$ for all $t\in[0,T)$.
It is worthy pointing out that, as mentioned in the Introduction, for the special case  Eq.(1.1) with $b=0$, one can apply the conserved quantity $\int_{\mathbb{R}} u m d x=\int_{\mathbb{R}} (u^2+u^2_x) d x$ and Sobolev's embedding theorem to uniformly bound $||u(t,\cdot)||_{L^\infty(\mathbb{R})}$. However, one cannot utilize any appropriate conservation laws of System (\ref{mn}) to control $||u(t,\cdot)||_{L^\infty}$ and $||v(t,\cdot)||_{L^\infty}$ directly. Anyway, the uniform boundedness for the solution $u$ to Eq.(1.1) with $b=0$ can be viewed as a special case of the above exponential increase assumption in finite time.

\begin{theorem4}
Suppose that $(m_0,n_0)\in H^s(\mathbb{R})\times H^s(\mathbb{R})$ with $s> \frac 1 2$
and $T>0$ be the maximal existence time of the corresponding solution $(m,n)$ to System (\ref{mn}) with $b=0$.
Assume that there exists a constant $C=C(||m_0||_{H^s},||n_0||_{H^s})>0$ such that
$||u(t,\cdot)||_{L^\infty},\,||v(t,\cdot)||_{L^\infty}\leq Ce^{Ct}$ for all $t\in[0,T)$.
Set $M(t)\triangleq M(t,q(t,y_0))=\inf\limits_{x\in\mathbb{R}} M(t,x)$ for some $y_0\in \mathbb{R}$, which is guaranteed by Remark 3.1 and Lemma 4.1, and $N(t)\triangleq (m+n)(t,q(t,y_0))$.
Let $m_0(x),\, n_0(x)\geq 0$ for all $x\in\mathbb{R}$,
and $m_0(x_0),\, n_0(x_0)> 0$ for some $x_0=q(0,y_0)$, which is ensured by Lemma 4.1.
If $M(0)<-2 C$ and
\begin{eqnarray}
\frac{M(0)}{N(0)}< -e^{(e^{C \eta}-1)}\left(\frac{2}{\eta N(0)}+1\right)+1,
\end{eqnarray}
where $\eta$ is the unique positive solution to the following equation w.r.t. $t$:
$$e^{(e^{C t}-1)}\left(\frac{C}{N(0)}e^{C t}+\frac{1}{2}C t e^{C t}+\frac{1}{2}\right)
+\frac{1}{2}\left(\frac{M(0)}{N(0)}-1\right)=0,
\quad t\geq 0.$$
Then the solution $(m,n)$ blows up at a time $T_0\in (0,\eta]$.
\end{theorem4}

\begin{proof}
In view of Remark 4.1, we here may assume $s\geq 3$ to prove the theorem. By (4.17), Lemma 4.3 and the assumption of the theorem, we have
\begin{eqnarray}
\frac{d}{dt} M(t)&=&\frac{d}{dt} M(t,q(t,y_0))\\ \nonumber
&=&(M_t+\frac{1}{2}(u v-u_x v_x)M_x)(t,q(t,y_0))\\ \nonumber
&\leq& -\frac{1}{2}M^2(t)+C e^{C t}N(t).
\end{eqnarray}
From System (\ref{mn}), we get
\begin{eqnarray}
\frac{d}{dt} N(t)&=&\frac{d}{dt} m(t,q(t,y_0))+\frac{d}{dt} n(t,q(t,y_0))\\ \nonumber
&=& -\frac{1}{2}M(t)N(t)+\frac{1}{2}(u v_x-u_x v)(m-n)(t,q(t,y_0)).
\end{eqnarray}
Note that (4.19)-(4.20) and the assumption imply $N(t)>0$ for all $t\in[0,T)$. Applying (4.28)-(4.29) and (4.31)-(4.32), one infers
\begin{eqnarray*}
&& N(t)\frac{d}{dt} M(t)-M(t)\frac{d}{dt} N(t)\\ \nonumber
&\leq& C e^{C t}N^2(t)-\frac{1}{2}(u v_x-u_x v)(m-n)(t,q(t,y_0))M(t)\\ \nonumber
&\leq& C e^{C t}N^2(t)+||u v||_{L^\infty}(||u||_{L^\infty}n+||v||_{L^\infty}m)N(t)\\ \nonumber
&\leq& C e^{C t}N^2(t),
\end{eqnarray*}
which gives
\begin{eqnarray*}
\frac{d}{dt} \left(\frac{M(t)}{N(t)}\right) \leq C e^{C t}.
\end{eqnarray*}
Integrating from $0$ to $t$ yields
\begin{eqnarray*}
\frac{M(t)}{N(t)}\leq \frac{M(0)}{N(0)}+ C\int_0^t  e^{C \tau}d\tau=\frac{M(0)}{N(0)}+e^{C t}-1,
\end{eqnarray*}
or hence
\begin{eqnarray}
M(t)\leq \left(\frac{M(0)}{N(0)}-1+e^{C t}\right)N(t).
\end{eqnarray}
On the other hand, according to (4.32) and Lemma 4.3, one gets
\begin{eqnarray*}
M(t)\geq \frac{-2}{N(t)}\frac{d}{dt}N(t)-C e^{C t},
\end{eqnarray*}
which along with (4.33) leads to
\begin{eqnarray}
\frac{d}{dt}\left(\frac{1}{N(t)}\right)\leq C e^{C t}\frac{1}{N(t)}
+\frac{1}{2}\left(\frac{M(0)}{N(0)}-1+e^{C t}\right).
\end{eqnarray}
Since $M(0)<-2 C$, $\frac{M(0)}{N(0)}<1$, if follows from Gronwall's inequality applied to (4.34) and the fact $x\leq e^{x}-1$ that
\begin{eqnarray}
0<\frac{1}{N(t)} \nonumber
&\leq& \frac{1}{N(0)} e^{\int_0^t C e^{C \tau}d\tau} \nonumber
+\frac{1}{2}\int_0^t e^{\int_{\tau}^t C e^{C s}d s}\left(\frac{M(0)}{N(0)}-1+e^{C \tau}\right)d\tau\\ \nonumber
&=&\frac{1}{N(0)} e^{(e^{C t}-1)}
+\frac{1}{2}\int_0^t e^{(e^{C t}-e^{C \tau})}e^{C \tau}d\tau\\ \nonumber
&&+\frac{1}{2}\int_0^t e^{\int_{\tau}^t C e^{C s}d s}\left(\frac{M(0)}{N(0)}-1\right)d\tau\\ \nonumber
&\leq& \frac{1}{N(0)} e^{(e^{C t}-1)}
+\frac{1}{2}\int_0^t e^{(e^{C t}-1)}d\tau+\frac{1}{2}\left(\frac{M(0)}{N(0)}-1\right)t\\ \nonumber
&=& e^{(e^{C t}-1)} \left(\frac{1}{N(0)}+\frac{1}{2}t\right)+\frac{1}{2}\left(\frac{M(0)}{N(0)}-1\right)t\\
&\triangleq& f(t).
\end{eqnarray}
Note that
$$f'(t)=e^{(e^{C t}-1)} \left(\frac{C}{N(0)}e^{C t}+\frac{1}{2}C t e^{C t}+\frac{1}{2}\right)
+\frac{1}{2}\left(\frac{M(0)}{N(0)}-1\right).$$
Since $M(0)<-2 C$ ensures $f'(0)=\frac{C}{N(0)}+\frac{M(0)}{2 N(0)}<0$, it follows from the facts
$$f{''}(t)=e^{(e^{C t}-1)}
\left(C e^{C t}(\frac{C}{N(0)}e^{C t}+\frac{C t}{2} e^{C t}+\frac{1}{2})
+C^2 e^{C t}(\frac{1}{N(0)}+\frac{t}{2})+\frac{C}{2} e^{C t}\right)>0$$
and $\lim\limits_{t\rightarrow +\infty} f'(t)=+\infty$ that there exists a unique $\eta>0$ such that $f'(\eta)=0$ and
$$f'(t)<0,\quad \text{if}\quad t\in [0,\eta) \quad \text{and}\quad f'(t)>0,\quad \text{if}\quad t>\eta.$$
By (4.30), we have $f(\eta)<0$. Noting that $f(0)=\frac{1}{N(0)}>0$ and $f(t)\in C[0,+\infty)$,
one can find a $T_0\in (0,\eta]$ such that
$$f(t)\rightarrow 0^{+},\quad \text{as}\,t\rightarrow T_0,$$
which together with (4.35) yields
$$N(t)\rightarrow +\infty,\quad \text{as}\,t\rightarrow T_0.$$
Thanks to (4.30) again, we get $\frac{M(0)}{N(0)}-1+e^{C\eta}<0$. This along with (4.33) ensures
$$\inf\limits_{x\in\mathbb{R}} (u_x n+v_x m)(t,x)=M(t)\rightarrow -\infty,\quad \text{as}\,t\rightarrow T_0.$$
According to Theorem 4.2, the solution $(m,n)$ blows up at the time $T_0\in(0,\eta]$. Therefore, we complete the proof of the theorem.
\end{proof}

\section{Peakon and weak kink solutions}

Remark 3.1 and the blow-up result Theorem 4.3 tell us that the solution to System (\ref{mn}) with sufficiently regular initial data may blow up in finite time. In this section, some explicit solutions to System (1.3) such as peakon and weak kink solutions are given, which demonstrates the occurrence of wave-breaking phenomena.
Let us first suppose that the single peakon solution of (1.3) with $b=0$ is of the following form
\begin{eqnarray}
u=c_1e^{-\mid x-ct\mid},\quad v=c_2e^{-\mid x-ct\mid}, \label{ocp}
\end{eqnarray}
where $c_1$ and $c_2$ are two constants.
Substituting (\ref{ocp}) into (1.3) with $b=0$ and integrating in the distribution sense,
we have
\begin{eqnarray}
c_1c_2=-3c. \label{C1}
\end{eqnarray}
In particular, for $c_1=c_2$, we recover the single peakon solution $u=\pm \sqrt{-3c}e^{-\mid x-ct\mid}$ of the cubic  CH equation \cite{Q1} with $b=0$ \cite{Gui-CMP,QXL}.

Next, let us show that System (1.3) with $b\neq0$ possesses a weak kink solution.
We assume System (1.3) has the following kink wave solution:
\begin{eqnarray}
u=C_1sgn(x-ct)\left(e^{-\mid x-ct\mid}-1\right),\quad
v=C_2sgn(x-ct)\left(e^{-\mid x-ct\mid}-1\right),
\label{kink1}
\end{eqnarray}
where the constant $c$ is the wave speed, and $C_1$ and $C_2$ are two constants.
In fact, if $C_1\neq0$ and $C_2\neq0$, the potentials $u$ and $v$ in (\ref{kink1}) are kink wave solutions due to
\begin{eqnarray}
\begin{array}{l}
\lim\limits_{x\rightarrow +\infty} u=-\lim\limits_{x\rightarrow -\infty} u=-C_1,\\
\lim\limits_{x\rightarrow +\infty} v=-\lim\limits_{x\rightarrow -\infty} v=-C_2.
\end{array}
\label{reasonkink}
\end{eqnarray}
We may easily check that the first order partial derivatives of (\ref{kink1}) read
\begin{eqnarray}
\begin{array}{l}
u_x=-C_1e^{-\mid x-ct\mid}, \quad u_t=c C_1e^{-\mid x-ct\mid},\\
v_x=-C_2e^{-\mid x-ct\mid}, \quad v_t=c C_2e^{-\mid x-ct\mid}.
\end{array}
\label{dkink1}
\end{eqnarray}
But, unfortunately, the second and higher order partial derivatives of (\ref{kink1}) do not exist
at the point $x=ct$.
Thus, like the case of peakon solutions, the kink wave solution in the form of (\ref{kink1}) should also be understood in the distribution sense, and therefor we call (\ref{kink1}) the weak kink solution to System (1.3) with $b\neq0$.

\begin{figure}
\centering
\includegraphics[width=2.2in]{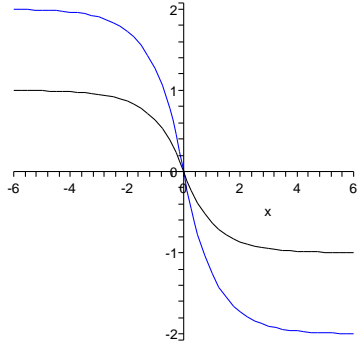}
\caption{\small{The weak kink solution at $t=0$. Black line: $u(x,0)$; Blue line: $v(x,0)$.}}
\label{figkink}
\end{figure}

Substituting (\ref{kink1}) and (\ref{dkink1}) into (1.3) and integrating through test functions,
we can arrive at
\begin{eqnarray}
\left\{\begin{array}{l}
c=-\frac{1}{2}b,\\
C_1C_2=-b.
\end{array}\right. \label{ikss}
\end{eqnarray}
In the above formula, $c=-\frac{1}{2}b$ means that the weak kink wave speed is exactly $-\frac{1}{2}b$.
In particular, we take $b=-2$ and $C_1=1$, then the corresponding weak kink solution is cast to
\begin{eqnarray}
u=sgn(x-t)\left(e^{-\mid x-t\mid}-1\right),\quad
v=2sgn(x-t)\left(e^{-\mid x-t\mid}-1\right).
\label{kink1s}
\end{eqnarray}
See Figure \ref{figkink} for the profile of the weak kink wave solution.

\bigskip
\noindent\textbf{Acknowledgments} \ Qiao was partially supported
by the Texas Norman Hackerman Advanced Research Program (Grant No.
003599-0001-2009) and supervised the U.S. Department of Education
UTPA-GAANN program (P200A120256). Yin was partially supported by
NNSFC (No. 11271382), RFDP (No. 20120171110014), and the key
project of Sun Yat-sen University (No. c1185).

\end{document}